\documentclass[A4paper,12pt]{book}
\usepackage{amssymb,amsmath}
\usepackage[russian]{babel}
\usepackage{graphicx}
\newcount\exnom\exnom=0

\long\def\exo/{\vspace{0.2cm} \noindent\advance\exnom by1{\bf
{\the\exnom}.}}

\newcommand{\ds}{\displaystyle}

\def\im{\mathop{\rm Im}\nolimits}

\begin{document}

\thispagestyle{empty}

\begin{center} \large \textbf{ENDOMORPHISM SEMIRINGS WITH ZERO\\ \vspace{1mm}
OF A FINITE SEMILATTICE OF A SPECIAL TYPE}
\end{center}

\begin{center} \textbf{Ivan Trendafilov}
\end{center}

\vspace{2mm}

\begin{quote}\centerline{{\bf Abstract}}
We investigate endomorphism semirings of a finite semilattice with one least element and one greatest element such that all the other elements form an antichain. We construct some new finite simple semirings.
\end{quote}

\textbf{Keywords:} endomorphism semiring, additively idempotent semiring, simple semiring, idempotent endomorphism, finite semilattice.

\textbf{MSC-2010:} 16Y60, 06A12, 20M10.

\vspace{5mm}

{\bf 1 \hspace{1mm} Introduction}

\vspace{3mm}

This paper is about the endomorphism semirings of a finite semilattice. There is a series of recent works where various problems of endomorphism semirings
have been considered, see [2], [5], [7] -- [12]. Here we consider  endomorphism semiring ${\mathcal{E}}_{\mathcal{\lozenge}_n}$ of the join-semilattice of $n$ -- element lattice with least element 0 and greatest element 1 and with a property that all the other elements form an antichain. We investigate zero-divisors, invertible elements and other elements of ${\mathcal{E}}_{\mathcal{\lozenge}_n}$ as well as some subsemirings and ideals. The semiring ${\mathcal{E}}_{\mathcal{\lozenge}_n}$ is not simple for $n \geq 5$ but there are many subsemirings which are simple. The study of simple algebras is a topic of great interest in algebra, moreover the study of finite simple semirings, has important application in cryptography.
In 2007, G. Maze, Ch. Monico and J. Rosenthal, see [6], first suggested some new ideas for
public key cryptography based on a concrete semigroup action built from simple semirings.
It is well known that  ``classical'' Pohlig-Hellman-type attacks
allows to compute the discrete logarithm in a cyclic group based on the Chinese remainder theorem. The complexity of the algorithm depends on the largest prime
factor of the order of this group. To prevent  Pohlig-Hellman-type attacks, the group can
be chosen to be a cyclic group of a large prime order,  that is the group is a
simple group. In order to prevent a similar attack on public key cryptosystems using semirings,
the use of simple semirings was suggested. If the used semiring was not be simple,
the semigroup action problem could be solved in a quotient semiring from which  an attacker
 may gain information to solve the semigroup action problem in the original
semiring.

The paper is organized as it follows. After the second section of preliminaries, in section 3 we consider the basic properties of semiring ${\mathcal{E}}_{\mathcal{\lozenge}_n}$, construct the addition and multiplication tables of semiring ${\mathcal{E}}_{\mathcal{\lozenge}_4}$ and prove that this semiring is simple. Here, we also construct two sorts of endomorphisms, defined by equalities (1) and (2) and prove that all these endomorphisms, the zero $\overline{0}$ and  $\overline{1}$ forms a subsemiring of ${\mathcal{E}}_{\mathcal{\lozenge}_n}$. In the next section we investigate the zero-divisors, the regular elements and the invertible elements of semiring ${\mathcal{E}}_{\mathcal{\lozenge}_n}$. Here, we extend the maximal ideal of the subsemiring of regular elements of ${\mathcal{E}}_{\mathcal{\lozenge}_n}$ to a maximal ideal of ${\mathcal{E}}_{\mathcal{\lozenge}_n}$. So, we prove that ${\mathcal{E}}_{\mathcal{\lozenge}_n}$ is not a simple semiring for $n > 4$. In section 5 we consider  endomorphisms $\alpha$ such that $\im(\alpha)$ is a subset of a fixed three-element chain $\mathcal{C}_i = \{0, a_i, 1\}$. We prove that for any $n \geq 4$ these endomorphisms form a simple semiring. Every semiring of this sort is a union of three semirings and one of them is also a simple semiring. In the next section we study the idempotent elements of semiring ${\mathcal{E}}_{\mathcal{\lozenge}_n}$. Since part of them, namely these endomorphisms which transform some of $a_i$ to $0$, has bad properties, we consider only idempotents having $1$ for a fixed point. We prove that these idempotents form a semiring and find some interesting subsemirings of him. By similar construction we obtain another semiring whose elements are also idempotents. In the last section we construct some simple semirings. The main results here are that  all the endomorphisms $\alpha$ such that $\im(\alpha) \subseteq \mathcal{C}_i\cup \mathcal{C}_j$ form a simple semiring (Theorem 7.1), but if $\ds \im(\alpha) \subseteq \bigcup_i^k \mathcal{C}_i$ and $k \geq 3$ the semiring is not simple (Theorem 7.5).

\vspace{5mm}

{\bf 2 \hspace{1mm} Preliminaries}

\vspace{3mm}

Let $\mathcal{M}$ be a \emph{semilattice (join-semilattice)} i.e. an algebra with binary operation $\vee$ such that

$\bullet$ $\; a\vee(b\vee c) = (a\vee b)\vee c$ for any $a, b, c \in \mathcal{M}$;
$\;\bullet$ $\; a\vee b = b \vee a$ for any $a, b\in \mathcal{M}$;

$\bullet$ $\; a\vee a = a$ for any $a \in \mathcal{M}$.

Another term used for arbitrary semilattice is a \emph{commutative idempotent semigroup} -- see [12].
For any $a, b\in \mathcal{M}$ we denote
$ a \leq b \; \iff \; a \vee b = b$.
In this notations, if there is a neutral element in the semilattice $\mathcal{M}$, it is the least element.

 Facts concerning semilattices can be found in [4].

 An algebra $R = (R,+,.)$ with two binary operations $+$ and $\cdot$ on $R$, is called a \emph{semiring} if

$\bullet\; (R,+)$ is a commutative semigroup,
$\bullet\; (R,\cdot)$ is a  semigroup,

$\bullet\;$ both distributive laws hold
$ x\cdot(y + z) = x\cdot y + x\cdot z$ and $(x + y)\cdot z = x\cdot z + y\cdot z$
for any $x, y, z \in R$.

Facts concerning semirings can be found in [3] and [7].
\vspace{1mm}

 Let $R = (R,+,.)$ be a semiring.

$\bullet\;$ If a neutral element $0$ of  semigroup $(R,+)$ exists and satisfies $0\cdot x = x\cdot 0 = 0$ for all $x \in R$, then it is called \emph{zero}.

$\bullet\;$ If a neutral element  of the semigroup $(R,\cdot)$ exists, it is called \emph{identity}.

For  semilattice $\mathcal{M}$ the set $\mathcal{E}_\mathcal{M}$ of the endomorphisms of $\mathcal{M}$ is a semiring
 with respect to the addition and multiplication defined with:

 $\bullet \; h = f + g \; \mbox{when} \; h(x) = f(x)\vee g(x) \; \mbox{for all} \; x \in \mathcal{M}$,

 $\bullet \; h = f\cdot g \; \mbox{when} \; h(x) = f\left(g(x)\right) \; \mbox{for all} \; x \in \mathcal{M}$.

 This semiring is called the \emph{ endomorphism semirimg} of $\mathcal{M}$.

An element $a$ of a semiring $R$ is called \emph{additively (multiplicatively) idempotent} if $a + a = a\;$ $\;(a\cdot a = a)$.
A semiring $R$ is called \emph{additively idempotent} if each of its elements is additively idempotent.

A semiring R with zero element $0$ is called \emph{zero-sum free} when for any $a, b \in R$ the equality  $a+b = 0$  implies
$a = b = 0$. Since every additively idempotent semiring is zero-sum free it follows that if $S$ is an additively idempotent semiring with zero $0$ then $S^* = S\backslash \{0\}$ is a subsemiring of $S$.

An element $a$ of a semiring $R$ is called \emph{additively (multiplicatively) absorbing element}  if and only if
$a + x = a\;\;$ $(a\cdot x = x\cdot a = a)\;\;$ for any $x \in R$.
The zero of $R$ is the unique multiplicative absorbing element; of course, it does not need to exist.
Following [7] an element of a semiring $R$
is called an \emph{infinity} if it is both additively and multiplicatively absorbing. Such an element we denote by $\infty$.

An element $x$ of semiring $R$ is called \emph{multiplicatively subidempotent} if and
only if $x^2 + x = x$, and a semiring $R$ is \emph{multiplicatively subidempotent semiring}
if and only if its elements are multiplicatively subidempotent.
Additively idempotent and  multiplicatively subidempotent semirings are called \emph{Viterbi semirings} and they play important
roles in modal logic, see [1].

 An equivalence relation $\sim$ on semiring $R$ is called \emph{congruence}
if it respects  semiring operations:
$x \sim y \;\;\; \mbox{implies} \;\;\; a + x \sim a + y$, $ax \sim ay$ and  $xa \sim ya$.
Semiring $R$ is called \emph{simple (congruence-simple)} if its only congruences are $\sim\; = id_R$
and $\sim \; = R\times R$.

Idempotent elements, nilpotent elements and zero-divisors in a semiring are definined in a similar way as  in the ring theory.

\vspace{6mm}

{\bf 3 \hspace{1mm} The endomorphism semiring with zero of a finite semilattice of a special type}

\vspace{3mm}

Let $\mathcal{\lozenge}_n$ be a finite lattice with  least element $0$, greatest element $1$ and all other elements $a_1, \ldots, a_{n-2}$ form an antichain (Fig. 1). We always assume that $n \geq 4$. So, the following equalities are fulfilled: ${a_i \vee a_i = a_i}$, $a_i \vee a_j = 1$, where $i \neq j$, $a_i \vee 0 = a_i$, $a_i \vee 1 = 1$, for any $i = 1, \ldots n-2$, $1 \vee 1 = 1$, $1 \vee 0 = 1$  and  $0 \vee 0 = 0$.

\begin{figure}[h]\centering
  \includegraphics[width=40mm]{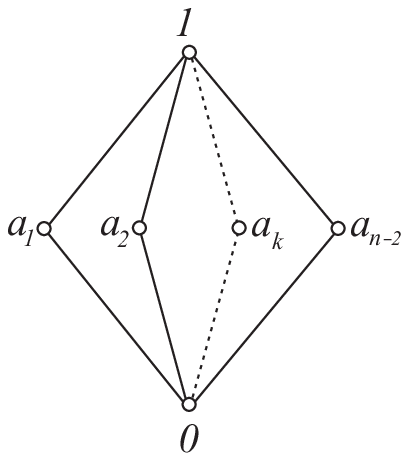}\\
\end{figure}
\vspace{-5mm}

\centerline{\small Figure 1. }

\vspace{1mm}

The endomorphism semiring of all endomorphisms of $\vee$ -- semilattice $\mathcal{\lozenge}_n$ having $0$ as a fixed point  is denoted by ${\mathcal{E}}_{\mathcal{\lozenge}_n}$.

Let the images of all the elements $a_i$ under the endomorphism $\alpha$ be elements $\alpha(a_i)\in \mathcal{\lozenge}_n$. Then $\alpha$ can be represented by ordered $n-1$ -- tuple $\wr\, \alpha(a_1), \alpha(a_2), \ldots \alpha(a_{n-2}), \alpha(1)\, \wr$.  So, the semiring  ${\mathcal{E}}_{\mathcal{\lozenge}_n}$ is additively idempotent with zero element $\wr\, 0, 0, \ldots, 0\, \wr$ and identity  $\wr{a_1, a_2 \ldots, a_{n-2}, 1}\wr$.

\vspace{3mm}

\textbf{Proposition 3.1} \textsl{Let $\alpha \in {\mathcal{E}}_{\mathcal{\lozenge}_n}$.}

\textsl{a.  The map $\alpha$ is isotone.}

\textsl{b. Let $\alpha \neq \wr\, 0, 0, \ldots, 0\, \wr$ and $\alpha(a_i) = 0$. Then for any $j \neq i$ it follows $\alpha(a_j) \neq 0$.
More precisely, $\alpha(a_i) = 0$ and $\alpha(1) = a_k$ imply $\alpha(a_j) = a_k$  for any $j \neq i$ and similarly $\alpha(a_i) = 0$ and $\alpha(1) = 1$ imply $\alpha(a_j) = 1$ for any $j \neq i$. }

\textsl{c. Let $\alpha(a_i) = a_k$ and $\alpha(1) = a_m$. Then either $\alpha(1) = a_k$ and $\alpha(a_j) = a_k$, for any $j \neq i$, or $\alpha(1) = a_k$, $\alpha(a_\ell) = 0$ for some $\ell$ and $\alpha(a_j) = a_k$ for every $j \neq i$ and $j \neq \ell$.}

\textsl{d. Let $\alpha(a_i) = a_k$ and $\alpha(1) = 1$. Then $\alpha(a_j) = 1$ for some elements $a_j$, where $j\neq i$ and $\alpha$ is a permutation of the remaining elements of some subset (maybe empty) of $\{a_1, \ldots, a_{n-2}\}$.}

\emph{Proof.} a. See p. 30 in [3].

b. Let for some $j \neq i$ it follows $\alpha(a_j) = 0$. Then $\alpha(1) = \alpha(a_i\vee a_j) = \alpha(a_i)\vee\alpha(a_j) = 0$ and this implies $\alpha = \wr\, 0, 0, \ldots, 0\, \wr$. Let $\alpha(a_i) = 0$ and $\alpha(1) = a_k$. Then $\alpha(a_j) = \alpha(a_i)\vee\alpha(a_j) = \alpha(1) = a_k$. Let $\alpha(a_i) = 0$ and $\alpha(1) = 1$. Then $\alpha(a_j) = \alpha(a_i)\vee\alpha(a_j) = \alpha(1) = 1$.

c. Let $\alpha(a_i) = a_k$. Then, see a. $\alpha(1) \neq 0$. If we assume $\alpha(1) = a_m$, then $a_m = \alpha(a_i)\vee\alpha(a_j) = a_k\vee\alpha(a_j)$. Now, there are two possibilities. The first one is $\alpha(a_j) = a_k = a_m$ for any $j \neq i$. The second one is: only for one $\ell$, using b., it follows $\alpha(a_\ell) = 0$ and $\alpha(a_j) = a_k = a_m$ for all $j \neq i$, $j \neq \ell$. Note that in all cases $\alpha(1) = a_k$.

d. Now $1 = \alpha(1) = \alpha(a_i)\vee\alpha(a_j) = a_k\vee\alpha(a_j)$. Hence, $\alpha(a_j) \neq a_k$ and $\alpha(a_j) \neq 0$. So, either $\alpha(a_j) = 1$ for some $a_j$, or $\alpha(a_j) = a_\ell$, where $\ell \neq k$. But in the last case $\alpha(a_r) = a_s$, where $r \neq i$ and $r \neq j$, it implies $s \neq k$ and $s \neq \ell$, and so on. Thus $\alpha$ is a permutation of some elements of the set $\{a_1, \ldots, a_{n-2}\}$.

\vspace{2mm}

Now we shall consider the ``$\,$least'' example of  semiring ${\mathcal{E}}_{\mathcal{\lozenge}_n}$.

\vspace{2mm}

\textbf{Example 3.2} Let $\mathcal{\lozenge}_4 = \; \left(\{0, a, b, 1\},\vee\right)$ is a four-element join-semilattice  having the folllowing $\vee$~table:
$\begin{array}{c|cccc}
\vee & 0 & a & b & 1\\ \hline
0 & 0 & a & b & 1\\
a & a & a & 1 & 1\\
b & b & 1 & b & 1\\
1 & 1 & 1 & 1 & 1\\
\end{array}$.
Let ${\mathcal{E}}_{\mathcal{\lozenge}_4}$ is an endomorphism semiring of this semilattice. From Proposition 3.1 it follows that ${\mathcal{E}}_{\mathcal{\lozenge}_4}$ has 16 elements and they are:
the zero element $\overline{0} =\wr\,0 0 0 \,\wr$, $\wr 0aa\wr$,  $\wr a0a \wr$, $\overline{a} = \wr aaa\wr$, $\wr 0bb\wr$, $\wr b0b\wr$, $\overline{b} = \wr bbb\wr$, the identity $i = \wr ab1\wr$, $\wr ba1\wr$, $\wr 011\wr$, $\wr 101\wr$, $\wr a11\wr$, $\wr 1a1\wr$, $\wr b11\wr$, $\wr 1b1\wr$ and  additively absorbing element $\overline{1} = \wr 111 \wr$.
 The addition and multiplication tables  are:
$$\footnotesize
    \begin{array}{c|cccccccccccccccc}
          + & \overline{0} & \wr 0aa \wr & \wr a0a \wr & \overline{a} & \wr 0bb\wr & \wr b0b\wr & \overline{b} & i &  \wr ba1 \wr & \wr 011 \wr & \wr 101\wr & \wr a11\wr & \wr 1a1\wr & \wr b11\wr  & \wr 1b1\wr  & \overline{1}\\ \hline
          \overline{0} & \overline{0} & \wr 0aa \wr & \wr a0a \wr & \overline{a} & \wr 0bb\wr & \wr b0b\wr & \overline{b} & i& \wr ba1 \wr & \wr 011 \wr & \wr 101\wr & \wr a11\wr & \wr 1a1\wr & \wr b11\wr  & \wr 1b1\wr  & \overline{1}\\
    \wr 0aa \wr & \wr 0aa \wr & \wr 0aa \wr & \overline{a} & \overline{a} & \wr 011\wr & \wr ba1\wr & \wr b11 \wr & \wr a11\wr &  \wr ba1 \wr & \wr 011 \wr & \wr 1a1\wr & \wr a11\wr & \wr 1a1\wr & \wr b11\wr  & \overline{1}  & \overline{1}\\
    \wr a0a \wr & \wr a0a \wr &  \overline{a}  & \wr a0a \wr & \overline{a} & i & \wr 101\wr & \wr 1b1 \wr & i  &  \wr 1a1 \wr & \wr a11 \wr & \wr 101\wr & \wr a11\wr & \wr 1a1\wr & \overline{1}  & \wr 1b1 \wr  & \overline{1}\\
    \overline{a}& \overline{a} &  \overline{a}  & \overline{a} & \overline{a} & \wr a11 \wr & \wr 1a1\wr & \overline{1} & \wr a11 \wr &  \wr 1a1 \wr & \wr a11 \wr & \wr 1a1\wr & \wr a11\wr & \wr 1a1\wr & \overline{1}  & \overline{1}  & \overline{1}\\
    \wr 0bb \wr & \wr 0bb \wr &  \wr 011 \wr & i & \wr a11 \wr & \wr 0bb \wr & \overline{b} & \overline{b} & i & \wr b11 \wr & \wr 011 \wr & \wr 1b1\wr & \wr a11\wr & \overline{1} & \wr b11 \wr  & \wr 1b1 \wr  & \overline{1}\\
    \wr b0b \wr & \wr b0b \wr &  \wr ba1 \wr & \wr 101 \wr & \wr 1a1 \wr & \overline{b} & \wr b0b \wr & \overline{b} & \wr 1b1\wr &  \wr ba1 \wr & \wr b11 \wr & \wr 101\wr & \overline{1} & \wr 1a1 \wr & \wr b11 \wr  & \wr 1b1 \wr  & \overline{1}\\
    \overline{b}& \overline{b} &  \wr b11 \wr & \wr 1b1 \wr & \overline{1} & \overline{b} & \overline{b} & \overline{b} & \wr 1b1 \wr &   \wr b11 \wr & \wr b11 \wr & \wr 1b1\wr & \overline{1} & \overline{1} & \wr b11 \wr  & \wr 1b1 \wr  & \overline{1}\\
   i & i &  \wr a11 \wr & i & \wr a11 \wr & i & \wr 1b1 \wr & \wr 1b1 \wr & i  &  \overline{1} & \wr a11 \wr & \wr 1b1\wr & \wr a11 \wr & \overline{1} & \overline{1}  & \wr 1b1 \wr  & \overline{1}\\
    \wr ba1 \wr & \wr ba1 \wr &  \wr ba1 \wr & \wr 1a1 \wr & \wr 1a1 \wr & \wr b11 \wr & \wr ba1 \wr & \wr b11 \wr & \overline{1} & \wr ba1 \wr & \wr b11 \wr & \wr 1a1 \wr & \overline{1} & \wr 1a1 \wr & \wr b11 \wr  & \overline{1}  & \overline{1}\\
   \wr 011 \wr & \wr 011 \wr &  \wr 011 \wr & \wr a11 \wr & \wr a11 \wr & \wr 011 \wr & \wr b11 \wr & \wr b11 \wr & \wr a11 \wr  &   \wr b11 \wr & \wr 011 \wr & \overline{1} & \wr a11 \wr & \overline{1} & \wr b11 \wr  & \overline{1}  & \overline{1}\\
   \wr 101 \wr & \wr 101 \wr &  \wr 1a1 \wr & \wr 101 \wr & \wr 1a1 \wr & \wr 1b1 \wr & \wr 101 \wr & \wr 1b1 \wr & \wr 1b1 \wr  &   \wr 1a1 \wr & \overline{1} & \wr 101 \wr & \overline{1} & \wr 1a1 \wr & \overline{1}  & \wr 1b1 \wr  & \overline{1}\\
   \wr a11 \wr & \wr a11 \wr &  \wr a11 \wr & \wr a11 \wr & \wr a11 \wr & \wr a11 \wr & \overline{1} & \overline{1} & \wr a11 \wr & \overline{1} & \wr a11 \wr & \overline{1} & \wr a11 \wr & \overline{1} & \overline{1}  & \overline{1}  & \overline{1}\\
   \wr 1a1 \wr & \wr 1a1 \wr &  \wr 1a1 \wr & \wr 1a1 \wr & \wr 1a1 \wr & \overline{1} & \wr 1a1 \wr & \overline{1} & \overline{1} &  \wr 1a1 \wr & \overline{1} & \wr 1a1 \wr & \overline{1} & \wr 1a1 \wr & \overline{1}  & \overline{1}  & \overline{1}\\
   \wr b11 \wr & \wr b11 \wr &  \wr b11 \wr & \overline{1} & \overline{1} & \wr b11 \wr & \wr b11 \wr & \wr b11 \wr & \overline{1} & \wr b11 \wr & \wr b11 \wr & \overline{1} & \overline{1} & \overline{1} & \wr b11 \wr  & \overline{1}  & \overline{1}\\
   \wr 1b1 \wr & \wr 1b1 \wr &  \overline{1} & \wr 1b1 \wr & \overline{1} & \wr 1b1 \wr & \wr 1b1 \wr & \wr 1b1 \wr & \wr 1b1 \wr  &  \overline{1} & \overline{1} & \wr 1b1 \wr & \overline{1} & \overline{1} & \overline{1}  & \wr 1b1 \wr & \overline{1}\\
    \overline{1} & \overline{1} &  \overline{1} & \overline{1} & \overline{1} & \overline{1} & \overline{1} & \overline{1}  & \overline{1} &  \overline{1} & \overline{1} & \overline{1} & \overline{1} & \overline{1} & \overline{1}  & \overline{1} & \overline{1}\\
   \end{array}$$

$$\footnotesize
    \begin{array}{c|cccccccccccccccc}
           \cdot & \overline{0} & \wr 0aa \wr & \wr a0a \wr & \overline{a} & \wr 0bb\wr & \wr b0b\wr & \overline{b} & i &  \wr ba1 \wr & \wr 011 \wr & \wr 101\wr & \wr a11\wr & \wr 1a1\wr & \wr b11\wr  & \wr 1b1\wr  & \overline{1} \\ \hline
    \overline{0}\vphantom{\int^0}&\overline{0} & \overline{0} & \overline{0} & \overline{0} & \overline{0} & \overline{0} & \overline{0} & \overline{0} & \overline{0}\vphantom{\int^0} & \overline{0} & \overline{0} & \overline{0} & \overline{0} & \overline{0}  & \overline{0}  & \overline{0}\\
    \wr 0aa \wr & \overline{0} & \overline{0} & \wr 0aa \wr & \wr 0aa \wr & \overline{0} & \wr 0bb\wr & \wr 0bb \wr & \wr 0aa \wr & \wr 0aa \wr & \overline{0}  & \wr 011\wr & \wr 0aa \wr & \wr 011\wr & \wr 0bb\wr  & \wr 011 \wr  & \wr 011 \wr\\
    \wr a0a \wr & \overline{0} &  \overline{0}  & \wr a0a \wr & \wr a0a \wr & \overline{0} & \wr b0b\wr & \wr b0b \wr & \wr a0a \wr &  \wr b0b \wr & \overline{0}& \wr 101\wr & \wr a0a\wr & \wr 101\wr & \wr b0b \wr  & \wr 101 \wr  & \wr 101 \wr\\
    \overline{a}& \overline{0} &  \overline{0}  & \overline{a} & \overline{a} & \overline{0} & \overline{b} & \overline{b} & \overline{a} &  \overline{b} & \overline{0} & \overline{1} & \overline{a} & \overline{1} & \overline{b}  & \overline{1}  & \overline{1}\\
    \wr 0bb \wr & \overline{0} &  \wr 0aa \wr & \overline{0} & \wr 0aa \wr & \wr 0bb \wr & \overline{0} & \wr 0bb \wr & \wr 0bb \wr & \wr 0aa \wr & \wr 011 \wr & \overline{0} & \wr 011\wr & \wr 0aa \wr & \wr 011 \wr  & \wr 0bb \wr  & \wr 011 \wr\\
    \wr b0b \wr & \overline{0} &  \wr a0a \wr & \overline{0} & \wr a0a \wr & \wr b0b \wr & \overline{0} & \wr b0b \wr & \wr b0b \wr  & \wr a0a \wr & \wr 101 \wr & \overline{0} & \wr 101 \wr & \wr a0a \wr & \wr 101 \wr  & \wr b0b \wr  & \wr 101 \wr \\
    \overline{b}& \overline{0} &  \overline{a} & \overline{0} & \overline{a} & \overline{b} & \overline{0} & \overline{b} & \overline{b} &  \overline{a} & \overline{1} & \overline{0} & \overline{1} & \overline{a} & \overline{1}  & \overline{b}  & \overline{1}\\
   i & \overline{0} & \wr 0aa \wr & \wr a0a \wr & \overline{a} & \wr 0bb\wr & \wr b0b\wr & \overline{b} & i  & \wr ba1 \wr & \wr 011 \wr & \wr 101\wr & \wr a11\wr & \wr 1a1\wr & \wr b11\wr  & \wr 1b1\wr  & \overline{1}\\
    \wr ba1 \wr & \overline{0} &  \wr a0a \wr & \wr 0aa \wr & \overline{a} & \wr b0b \wr & \wr 0bb \wr & \overline{b} & \wr ba1 \wr &  i & \wr 101 \wr & \wr 011 \wr & \wr 1a1 & \wr a11 \wr & \wr 1b1 \wr  & \wr b11 \wr  & \overline{1}\\
   \wr 011 \wr & \overline{0} &  \wr 0aa \wr & \wr 0aa \wr & \wr 0aa \wr & \wr 0bb \wr & \wr 0bb \wr & \wr 011 \wr & \wr 011 \wr  &  \wr 011 \wr & \wr 011 \wr & \wr 011 \wr & \wr 011 \wr & \wr 011 \wr & \wr 011 \wr  & \wr 011 \wr  & \wr 011 \wr\\
   \wr 101 \wr & \overline{0} &  \wr a0a \wr & \wr a0a \wr & \wr a0a \wr & \wr b0b \wr & \wr b0b \wr & \wr b0b \wr & \wr 101 \wr & \wr 101 \wr & \wr 101 \wr & \wr 101 \wr & \wr 101 \wr  & \wr 101 \wr & \wr 101 \wr  & \wr 101 \wr  & \wr 101 \wr\\
   \wr a11 \wr & \overline{0} &  \wr 0aa \wr & \overline{a} & \overline{a} & \wr 0bb \wr & \overline{b} & \overline{b} & \wr a11 \wr &  \wr b11 \wr & \wr 011 \wr & \overline{1} & \wr a11 \wr & \overline{1} & \wr b11 \wr  & \overline{1}  & \overline{1}\\
   \wr 1a1 \wr & \overline{0} &  \wr a0a \wr & \overline{a} & \overline{a} & \wr b0b \wr & \overline{b} & \overline{b} & \wr 1a1 \wr &  \wr 1b1 \wr & \wr 101 \wr & \overline{1} & \wr 1a1 \wr & \overline{1} & \wr 1b1 \wr  & \overline{1}  & \overline{1}\\
   \wr b11 \wr & \overline{0} &  \overline{a} & \wr 0aa \wr & \overline{a} & \overline{b} & \wr 0bb \wr & \overline{b} & \wr b11 \wr  &  \wr a11 \wr & \overline{1} & \wr 011 \wr & \overline{1} & \wr a11 \wr & \overline{1}  & \wr b11 \wr  & \overline{1}\\
   \wr 1b1 \wr & \overline{0} &  \overline{a} & \wr a0a \wr & \overline{a} & \overline{b} & \wr b0b \wr & \overline{b} & \wr 1b1 \wr &  \wr 1a1 \wr & \overline{1} & \wr 101 \wr & \overline{1} & \wr 1a1 \wr & \overline{1}  & \wr 1b1 \wr & \overline{1}\\
    \overline{1} & \overline{0} &  \overline{a} & \overline{a} & \overline{a} & \overline{b} & \overline{b} & \overline{b}  & \overline{1}  & \overline{1} & \overline{1} & \overline{1} & \overline{1} & \overline{1} & \overline{1}  & \overline{1} & \overline{1}\\
   \end{array}$$

\vspace{3mm}

Now we prove that there are not any proper ideals in semiring ${\mathcal{E}}_{\mathcal{\lozenge}_4}$, i.e. it is a simple semiring.
Let $I$ be an ideal of  semiring ${\mathcal{E}}^{[\,i\,]}_{\mathcal{\lozenge}_4}$. Since $\overline{0}$ is a zero element of this semiring, we may suppose that $\overline{0} \in I$. We shall consider the following four cases:

\vspace{1mm}

\emph{Case 1.} Let us assume that $\wr 011 \wr \in I$. Then $\wr 011 \wr\cdot \wr b0b \wr = \wr 0bb \wr \in I$, $\wr b0b \wr\cdot \wr 011 \wr = \wr 101 \wr\in I$ and $\wr 101 \wr\cdot \wr 0aa \wr = \wr a0a \wr\in I$. But $\wr a0a \wr + \wr 0bb \wr = i \in I$ that is $I = {\mathcal{E}}_{\mathcal{\lozenge}_4}$.

\vspace{1mm}

\emph{Case 2.} Let us assume that $\alpha \in I$, where $\alpha(1) = 1$ and $\alpha \neq i$. Then $\wr 011 \wr\cdot \alpha = \wr 011 \wr \in I$ and we go to Case 1.

\vspace{1mm}

\emph{Case 3.} Let us assume that $\alpha \in I$, where $\alpha(1) = a$. Then $\alpha\cdot \wr 101 \wr = \beta \in I$, where $\beta(1) = 1$, and we go to Case 2.

\vspace{1mm}

\emph{Case 4.} Let us assume that $\alpha \in I$, where $\alpha(1) = b$. Then $\alpha\cdot \wr 011 \wr = \beta \in I$, where $\beta(1) = 1$, and we go to Case 2.

\vspace{1mm}

Thus, we find that either $I = \{\overline{0}\}$, or $I = {\mathcal{E}}_{\mathcal{\lozenge}_4}$.

\vspace{3mm}

 Let $\alpha_i \in {\mathcal{E}}_{\mathcal{\lozenge}_n}$, where $\varphi_i(a_j) = a_i$ and $\varphi_i(1) = a_i$ for all $i, j = 1, \ldots, n - 2$. As in Example 3.2 we denote $\alpha_i = \overline{a_i}$. By the same way $\overline{0} = \wr 0, 0, \ldots, 0 \wr$ and $\overline{1} = \wr 1, 1, \ldots, 1 \wr$.  Endomorphisms $\overline{a_i}$, $i = 1, \ldots, n - 2$, $\overline{0}$ è $\overline{1}$ are called \textbf{\emph{almost constant}}. The set of all the almost constant endomorphisms is denoted by $\mathcal{AC}\left({\mathcal{E}}_{\mathcal{\lozenge}_n}\right)$.

\vspace{3mm}

From equalities:
$$ \overline{0} + \overline{a_i} = \overline{a_i},\; \overline{0} + \overline{1} = \overline{1},\; \overline{a_i} + \overline{1} = \overline{1},\; \overline{a_i} + \overline{a_i} = \overline{a_i},\; \overline{a_i} + \overline{a_j} = \overline{1}, \; \mbox{where}\; i \neq j,$$
$$\overline{0}\, \cdot \,\wr\, \alpha(a_1), \ldots, \alpha(1) \,\wr = \wr\, \alpha(a_1), \ldots, \alpha(1) \,\wr \, \cdot \, \overline{0} = \overline{0},    $$
$$\overline{a_i}\, \cdot \,\wr\, \alpha(a_1), \ldots, \alpha(a_i), \ldots \alpha(1) \,\wr = \wr\, \alpha(a_i), \ldots, \alpha(a_i) \,\wr = \overline{\alpha(a_i)},$$
$$\overline{1}\, \cdot \,\wr\, \alpha(a_1), \ldots, \alpha(a_i), \ldots \alpha(1) \,\wr = \wr\, \alpha(1), \ldots, \alpha(1) \,\wr = \overline{\alpha(1)}.$$
it follows

\vspace{3mm}

\textbf{Proposition 3.3} \textsl{For any $n \geq 4$ the set $\mathcal{AC}\left({\mathcal{E}}_{\mathcal{\lozenge}_n}\right)$ is a right ideal of  semiring ${\mathcal{E}}_{\mathcal{\lozenge}_n}$ .}

\vspace{3mm}

Note that all the elements of $\mathcal{AC}\left({\mathcal{E}}_{\mathcal{\lozenge}_n}\right)$ are obviously idempotents.

\vspace{3mm}

\textbf{Remark 3.4}  Let us compare the two semirings ${\mathcal{E}}_{\mathcal{\lozenge}_n}$ (with zero) and $\widehat{\mathcal{E}}_{\mathcal{A}_n}$ (without zero) from the paper [10]. It is easy to see that there is a bijection between  set $\mathcal{AC}\left({\mathcal{E}}_{\mathcal{\lozenge}_n}\right)$ of almost constant elements of ${\mathcal{E}}_{\mathcal{\lozenge}_n}$ and the set of constant elements of $\widehat{\mathcal{E}}_{\mathcal{A}_n}$. But the last set is an ideal of the semiring $\widehat{\mathcal{E}}_{\mathcal{A}_n}$ while $\mathcal{AC}\left({\mathcal{E}}_{\mathcal{\lozenge}_n}\right)$ is only a right ideal of semiring ${\mathcal{E}}_{\mathcal{\lozenge}_n}$. Moreover, the set ${\mathcal{E}}_{\mathcal{A}_n}$ of all endomorphisms which are non constant is a subsemiring of $\widehat{\mathcal{E}}_{\mathcal{A}_n}$ while the set of non almost constant endomorphisms of ${\mathcal{E}}_{\mathcal{\lozenge}_n}$ is not closed under the addition as well as under the multiplication, see Example 3.2.

\vspace{3mm}

Let us consider  endomorphisms  $\varphi_i \in {\mathcal{E}}_{\mathcal{\lozenge}_n}$, such that
$$\varphi_{i}(a_j) = \left\{ \begin{array}{l} 0, \; \mbox{for}\; i = j\\
1, \; \mbox{for}\; i \neq j \end{array} \right., \eqno{(1)}$$
where $i, j \in \{1, \ldots, n-2\}$. Obviously $\varphi_i(1) = 1$.

It is easy to see that $\varphi_i + \varphi_i = \varphi_i$ and $\varphi_i + \varphi_j = \overline{1}$ for any $i, j \in \{1, \ldots, n-2\}$ and $i \neq j$. Similarly, it follows $\varphi_i \cdot \varphi_i = \varphi_i$ and $\varphi_i \cdot \varphi_j = \varphi_i$ for any $i, j \in \{1, \ldots, n-2\}$ and $i \neq j$. Let us denote by ${\mathcal{E}}^{[\,0, 1\,]}_{\mathcal{\lozenge}_n}$ the set of all endomorphisms $\varphi_i$, $i = 1, \ldots, n-2$ and endomorphisms $\overline{0}$ and $\overline{1}$. Obviously $\varphi_i + \overline{0} = \varphi_i$, $\;\varphi_i\cdot \overline{0} = \overline{0}\cdot \varphi_i = \overline{0}$. Since  $\varphi_i + \overline{1} = \varphi_i \cdot \overline{1}  = \overline{1}\cdot \varphi_i = \overline{1}$, it follows that in semiring  $\left({\mathcal{E}}^{[0,1]}_{\mathcal{\lozenge}_n}\right)^* = {\mathcal{E}}^{[0,1]}_{\mathcal{\lozenge}_n}\backslash \{\overline{0}\}$   we have $\overline{1} = \infty$. Thus we prove

\vspace{3mm}

\textbf{Lemma 3.5} \textsl{For any $n \geq 4$ the set ${\mathcal{E}}^{[\,0, 1\,]}_{\mathcal{\lozenge}_n}$ is a subsemiring  of the semiring ${\mathcal{E}}_{\mathcal{\lozenge}_n}$ with a zero element $\overline{0}$. The element $\overline{1}$ is an infinity of $\left({\mathcal{E}}^{[\,0, 1\,]}_{\mathcal{\lozenge}_n}\right)^*$.}

\vspace{3mm}

Since  all the elements of  semiring ${\mathcal{E}}^{[\,0, 1\,]}_{\mathcal{\lozenge}_n}$ are both additively and multiplicatively idempotents it immediately follows

\vspace{3mm}

\textbf{Corollary 3.6} \textsl{The semiring ${\mathcal{E}}^{[\,0, 1\,]}_{\mathcal{\lozenge}_n}$ is a Viterbi semiring.}

\vspace{3mm}

Let us also consider the endomorphisms $\psi_{i,j} \in {\mathcal{E}}_{\mathcal{\lozenge}_n}$, such that
$$\psi_{i,j}(a_k) = \left\{ \begin{array}{l} a_j, \; \mbox{for}\; i = k\\
1, \; \mbox{for}\; i \neq k \end{array},\right. \eqno{(2)}$$
where $i, j, k \in \{1, \ldots, n-2\}$. Obviously, $\psi_{i,j}(1) = 1$.

Let us denote by ${\mathcal{E}}^{[\,a, 1\,]}_{\mathcal{\lozenge}_n}$ the set of all endomorphisms $\psi_{i,j}$, $i, j = 1, \ldots, n-2$ and endomorphisms $\overline{0}$ and $\overline{1}$.

\vspace{3mm}

\textbf{Lemma 3.7} \textsl{The set ${\mathcal{E}}^{[\,a, 1\,]}_{\mathcal{\lozenge}_n}$ is a subsemiring  of  semiring ${\mathcal{E}}_{\mathcal{\lozenge}_n}$ with a zero element $\overline{0}$.  The element $\overline{1}$ is an infinity of $\left({\mathcal{E}}^{[\,a, 1\,]}_{\mathcal{\lozenge}_n}\right)^*$.}

\emph{Proof.} Now we  calculate
$$\psi_{i,j} + \psi_{i,j} = \psi_{i,j}, \; \psi_{i,j} + \psi_{k,\ell} = \overline{1}, \; \mbox{where either}\; k \neq i, \; \mbox{or}\; k = i\; \mbox{and}\; j \neq \ell,$$
$i, j, k, \ell \in \{1, \ldots, n-2\}$ and also $\psi_{i,j} + \overline{0} = \psi_{i,j}$, $\,\psi_{i,j} + \overline{1} = \overline{1}$ and $\overline{0} + \overline{1} = \overline{1}$.

For arbitrary $a_m$ it follows
$$(\psi_{i,j}\cdot\psi_{j,\ell})(a_m) = \psi_{j,\ell}(\psi_{i,j}(a_m)) = \left\{ \begin{array}{ll} \psi_{j,\ell}(a_j) = a_\ell, &  \mbox{for}\; m = i\\
\psi_{j,\ell}(1) = 1, & \mbox{for}\; m \neq i \end{array}\right. = \psi_{i,\ell}(a_m).$$
Therefore, $\psi_{i,j}\cdot\psi_{j,k} = \psi_{i,k}$. Similarly $\psi_{i,j}\cdot\psi_{k,\ell} = \overline{1}$ for $j \neq k$ and $\psi_{i,j} \cdot \overline{1} = \overline{1}\cdot \psi_{ij} = \overline{1}$. Thus we prove that $\overline{1} = \infty$ in semiring $\left({\mathcal{E}}^{[\,a, 1\,]}_{\mathcal{\lozenge}_n}\right)^*$.

\vspace{3mm}

 All the endomorphisms $\varphi_i$, $i = 1, \ldots, n-2$ defined by (1) and $\psi_{i,j}$, $i, j = 1, \ldots, n - 2$ defined by (2), are called \textbf{\emph{almost absorbing}}. The subset of  semiring $\mathcal{E}_{\mathcal{\lozenge}_n}$ consisting of all the almost absorbing endomorphisms and endomorphisms $\overline{0}$ and $\overline{1}$, is denoted by $\mathcal{AA}\left({\mathcal{E}}_{\mathcal{\lozenge}_n}\right)$.

\vspace{3mm}

\textbf{Theorem 3.8} \textsl{For any $n \geq 4$ the set $\mathcal{AA}\left({\mathcal{E}}_{\mathcal{\lozenge}_n}\right)$ is a subsemiring of  semiring ${\mathcal{E}}_{\mathcal{\lozenge}_n}$ with a zero element $\overline{0}$. The element $\overline{1}$ is an infinity of semiring  $\left(\mathcal{AA}\left({\mathcal{E}}_{\mathcal{\lozenge}_n}\right)\right)^*$. The semiring ${\mathcal{E}}^{[\,0, 1\,]}_{\mathcal{\lozenge}_n}$ is an ideal of semiring $\mathcal{AA}\left({\mathcal{E}}_{\mathcal{\lozenge}_n}\right)$.}

\emph{Proof.} We observe that $\varphi_k + \psi_{i,j} = \psi_{i,j}$ for any $i, j, k \in \{1, \ldots, n-2\}$.     After using lemmas 3.5 and 3.7, we calculate
$$\varphi_k \cdot \psi_{i,j} = \varphi_k, \; \psi_{i,j} \cdot \varphi_j = \varphi_j, \; \psi_{i,j} \cdot \varphi_k = \overline{1} \; \mbox{for}\;\; k \neq j,$$
where $i, j, k \in \{1, \ldots, n-2\}$, and this completes the proof.

\vspace{3mm}

\textbf{Remark 3.9} The semiring $\mathcal{AA}\left({\mathcal{E}}_{\mathcal{\lozenge}_n}\right)$ is not an ideal of  ${\mathcal{E}}_{\mathcal{\lozenge}_n}$. For instance
$$ \wr\,0, a_1, \ldots a_1\, \wr \, \cdot \, \wr \, a_1, 1, \ldots, 1\, \wr = \wr \, a_1, 1, \ldots, 1\, \wr \, \cdot \, \wr\,0, a_1, \ldots a_1\, \wr = \wr\,0, a_1, \ldots a_1\, \wr.$$

Let us consider the subset of $\mathcal{AA}\left({\mathcal{E}}_{\mathcal{\lozenge}_n}\right)$ consisting of all the endomorphisms $\varphi_i$ and $\psi_{i,i}$, where $i \in \{1, \ldots, n-2\}$ and $\overline{1}$. We denote this set by $\mathcal{ID}\left(\mathcal{AA}\left({\mathcal{E}}_{\mathcal{\lozenge}_n}\right)\right)$. Since $(\psi_{i,j})^2 = \overline{1}$ for $i \neq j$, $\,\mathcal{ID}\left(\mathcal{AA}\left({\mathcal{E}}_{\mathcal{\lozenge}_n}\right)\right)$  is the set of all idempotent elements of semiring $\mathcal{AA}\left({\mathcal{E}}_{\mathcal{\lozenge}_n}\right)$.

\vspace{3mm}

\textbf{Proposition 3.10} \textsl{The set $\mathcal{ID}\left(\mathcal{AA}\left({\mathcal{E}}_{\mathcal{\lozenge}_n}\right)\right)$ is a subsemiring of semiring $\mathcal{AA}\left({\mathcal{E}}_{\mathcal{\lozenge}_n}\right)$.}

\emph{Proof.} From the proof of Theorem 3.8 it follows $\varphi_k + \psi_{i,i} = \psi_{i,i}$, $\varphi_k\cdot \psi_{i,i} = \varphi_k$, $\psi_{i,i} \cdot \varphi_k = \overline{1}$ for $i \neq k$ and $\psi_{i,i} \cdot \varphi_i = \varphi_i$, where $i, k \in \{1, \ldots, n-2\}$. From the proof of Lemma 3.7 it follows $\psi_{i,i}\cdot \psi_{i,i} = \psi_{i,i}$ and $\psi_{i,i}\cdot \psi_{k,k} = \psi_{k,k}\cdot \psi_{i,i} = \overline{1}$, where $i, k \in \{1, \ldots, n-2\}$, $i \neq k$.

\vspace{3mm}

Note that $\mathcal{ID}\left(\mathcal{AA}\left({\mathcal{E}}_{\mathcal{\lozenge}_n}\right)\right)$ is not an ideal of $\mathcal{AA}\left({\mathcal{E}}_{\mathcal{\lozenge}_n}\right)$ because $\psi_{i,i}\cdot \psi_{i,j} = \psi_{i,j}$ for $i \neq j$, where $i, j \in \{1, \ldots, n-2\}$.

\vspace{7mm}

{\bf 4 \hspace{1mm} Zero-divisors and invertible elements of the endomorphism semiring}

\vspace{3mm}

Let us consider  endomorphisms $\alpha_{\,0,j}^{(i)}$ such that $$\alpha_{\,0,j}^{(i)}(a_k) = \left\{\begin{array}{ll} 0, & \mbox{for}\; k = j\\ a_i, & \mbox{for}\; k \neq j \end{array}\right., \eqno{(3)}$$
where $i, j, k \in \{1, \ldots, n-2\}$.
\vspace{3mm}

To find nilpotent elements in semiring ${\mathcal{E}}_{\mathcal{\lozenge}_n}$ we observe that, when for some $\alpha \in {\mathcal{E}}_{\mathcal{\lozenge}_n}$ we have $1 \in \im(\alpha)$, using Proposition 3.1 a., it follows $\alpha(1) = 1$, that is $1$ is a fixed point of $\alpha$ and $\alpha$ is not a nilpotent element. So, for nilpotent endomorphism $\alpha$ follows $\alpha(1) = a_i$ for some $i = 1, \ldots, n-2$. From Proposition 3.1 c. we conclude that either $\alpha = \overline{a_i}$, but $\overline{a_i}$ is an idempotent, or $\alpha(a_k) = 0$ for some $k = 1, \ldots, n-2$. Note, that  endomorphism $\alpha$ such that $\alpha(a_j) = a_i$ for $j \neq k$ and $\alpha(a_k) = 0$ is also idempotent in all cases when $k \neq i$. So, only  endomorphism $\alpha_{0,i} = \alpha_{\,0,i}^{(i)}$, see (3), such that
$$\alpha_{0,i}(a_k) = \left\{\begin{array}{ll} 0, & \mbox{for}\; k = i\\ a_i, & \mbox{for}\; k \neq i \end{array}\right. \eqno{(4)}$$
is a nilpotent element. Hence, all nilpotent endomorphisms of ${\mathcal{E}}_{\mathcal{\lozenge}_n}$ are
$$\wr\, 0, a_1, \ldots, a_1\, \wr, \wr\, a_2, 0, a_2, \ldots, a_2\, \wr, \ldots, \wr\, a_{n-2}, \ldots, a_{n-2}, 0, a_{n-2}\,\wr.$$

Note that the sum and product of any two of them is not a nilpotent element.

Since $\wr\, a_1, 0, a_1, \ldots, a_1\,\wr\cdot \wr\, 0, 1, \ldots, 1\, \wr = \overline{0}$ it follows that there are zero-divisors of ${\mathcal{E}}_{\mathcal{\lozenge}_n}$ which are not nilpotent elements.

Let us note that $\overline{1}$ is not a zero-divisor, but all other almost constant endomorphisms, different from $\overline{0}$, are zero-divisors, for instance $\overline{a_i}\cdot \alpha_{0,i} = \overline{0}$, where $i = 1, \ldots, n-2$, see (4).

\vspace{3mm}

\textbf{Proposition 4.1} \textsl{Any endomorphism $\alpha \in {\mathcal{E}}_{\mathcal{\lozenge}_n}$, which is not an almost constant endomorphism, is a zero-divisor if and only if $0 \in \im(\alpha)$.}

\emph{Proof.} We assume that $0 \notin \im(\alpha)$. Then, if $\alpha(1) = a_i$, from Proposition 3.1 c. it follows that $\alpha = \overline{a_i}$. So, $\alpha(1) = 1$ and for arbitrary $\beta \in {\mathcal{E}}_{\mathcal{\lozenge}_n}$, $\beta \neq \overline{0}$ we have $(\alpha\cdot \beta)(1) = \beta(\alpha(1)) = \beta(1) \neq 0$, i.e. $\alpha$ is not a left zero-divisor. Let $\beta \in {\mathcal{E}}_{\mathcal{\lozenge}_n}$, $\beta \neq \overline{0}$ and $\beta(a_i) \neq 0$, where $i = 1, \ldots, n-2$. Then either $\beta(a_i) = 1$ and $(\beta\cdot \alpha)(a_i) = \alpha(\beta(a_i)) = 1$, or $\beta(a_i) = a_k$ and $(\beta\cdot \alpha)(a_i) = \alpha(\beta(a_i)) = \alpha(a_k) \neq 0$. So, $\alpha$ is not a right zero-divisor.

Conversely, we assume that $0 \in \im(\alpha)$. If $1 \in \im(\alpha)$, from Proposition 3.1 b. we have $\alpha = \varphi_i$ for some $i = 1, \ldots, n-2$ and then $\alpha_{\,0j}^{(i)}\cdot \varphi_i = \overline{0}$, where $j = 1, \ldots, n-2$. If $1 \notin \im(\alpha)$, then $\alpha = \alpha_{\,0j}^{(i)}$ for some $i$ and then $\alpha_{\,0j}^{(i)}\cdot \alpha_{0i} = \overline{0}$, where $j = 1, \ldots, n-2$.

\vspace{3mm}

An element of ${\mathcal{E}}_{\mathcal{\lozenge}_n}$ that is neither $0$, nor a zero-divisor is called regular. From Proposition 4.1 it follows that the regular elements of the semiring are  endomorphisms $\alpha \in {\mathcal{E}}_{\mathcal{\lozenge}_n}$ such that $0 \notin \im(\alpha)$ and also $\alpha \neq \overline{a_i}$, where $i = 1, \ldots, n-2$. Let us denote by $\mathfrak{Reg}\left(\mathcal{E}_{\mathcal{\lozenge}_n}\right)$ the set of regular endomorphisms of ${\mathcal{E}}_{\mathcal{\lozenge}_n}$.

\vspace{3mm}

\textbf{Proposition 4.2} \textsl{For any $n \geq 4$ the set $\mathfrak{Reg}\left(\mathcal{E}_{\mathcal{\lozenge}_n}\right)$ is a subsemiring of ${\mathcal{E}}_{\mathcal{\lozenge}_n}$ without zero. The identity $i$ is both an additively neutral and  a multiplicatively neutral element
 of $\mathfrak{Reg}\left(\mathcal{E}_{\mathcal{\lozenge}_n}\right)$ and the endomorphism $\overline{1}$ is infinity.}

\emph{Proof.} From Proposition 3.1 b., c. and d. it follows that for every $\alpha \in \mathfrak{Reg}\left(\mathcal{E}_{\mathcal{\lozenge}_n}\right)$ we obtain $\alpha(1) = 1$ and for any $a_i$ either $\alpha(a_i) = a_j$, or $\alpha(a_i) = 1$ where $i,j = 1, \ldots, n-2$. Now let $\beta \in \mathfrak{Reg}\left(\mathcal{E}_{\mathcal{\lozenge}_n}\right)$. Then $(\alpha + \beta)(1) = \alpha(1) \vee \beta(1) = 1$, $(\alpha\cdot \beta)(1) = \beta(\alpha(1)) = 1$ and $(\alpha + \beta)(a_i) = \alpha(a_i) \vee \beta(a_i) = 1$. Since $(\alpha\cdot \beta)(a_i) = \beta(\alpha(a_i)) = \beta(a_j)$, it follows that $(\alpha\cdot \beta)(a_i) = a_k$ or $(\alpha\cdot \beta)(a_i) = 1$. Hence, we prove that $\alpha + \beta, \alpha\cdot \beta \in \mathfrak{Reg}\left(\mathcal{E}_{\mathcal{\lozenge}_n}\right)$. So, $\mathfrak{Reg}\left(\mathcal{E}_{\mathcal{\lozenge}_n}\right)$ is a subsemiring of ${\mathcal{E}}_{\mathcal{\lozenge}_n}$. It is clear that identity $i$ is the least element $\mathfrak{Reg}\left(\mathcal{E}_{\mathcal{\lozenge}_n}\right)$ which means that $i$ is a neutral element of  additive semigroup $\left(\mathfrak{Reg}\left(\mathcal{E}_{\mathcal{\lozenge}_n}\right),+\right)$. Consequently there is not zero in the semiring $\mathfrak{Reg}\left(\mathcal{E}_{\mathcal{\lozenge}_n}\right)$. Since $\alpha(1) = 1$ implies $\alpha\cdot \overline{1} = \overline{1}\cdot \alpha = \overline{1}$ and also $\alpha + \overline{1} = \overline{1}$, for any $\alpha \in \mathfrak{Reg}\left(\mathcal{E}_{\mathcal{\lozenge}_n}\right)$ we have  $\overline{1} = \infty$.

\vspace{3mm}

We now seek the invertible elements in semiring ${\mathcal{E}}_{\mathcal{\lozenge}_n}$ which are  invertible regular endomorphisms. From Proposition 3.1 d.
 follows that the endomorphisms $\alpha = \wr\, a_{k_1},\,a_{k_2},\, \ldots, a_{k_{n-2}},\, 1 \, \wr$, where  $k_s \in \{1, \ldots, n - 2\}$ for any $s = 1, \ldots, n - 2$, are all the permutations of  elements $a_1, \ldots, a_{n-2}$. It is clear that the set of these permutations is a subgroup of the semigroup $\left(\mathfrak{Reg}\left(\mathcal{E}_{\mathcal{\lozenge}_n}\right),\cdot\right)$ and this group is isomorphic of the symmetric group $S_{n-2}$. We denote this group by $\mathcal{P}\left(\mathcal{E}_{\mathcal{\lozenge}_n}\right)$.

Let $\alpha$ be an invertible element of the semiring $\mathfrak{Reg}\left(\mathcal{E}_{\mathcal{\lozenge}_n}\right)$. Obviously, $\alpha(a_i) \neq 1$ for any $i = 1, \ldots, n-2$. So, $\alpha$ is a permutation on the set $\{a_1, \ldots, a_n\}$, that is the set of the invertible elements of $\mathfrak{Reg}\left(\mathcal{E}_{\mathcal{\lozenge}_n}\right)$ is  group $\mathcal{P}\left(\mathcal{E}_{\mathcal{\lozenge}_n}\right)$.

\vspace{3mm}

Let us consider  regular endomorphisms $\alpha$ that satisfy  condition $\alpha(a_i) = 1$ for some $i = 1, \ldots, {n-2}$ and denote the set of such endomorphisms by $\mathcal{M}\left(\mathfrak{Reg}\left(\mathcal{E}_{\mathcal{\lozenge}_n}\right)\right)$. So, semiring $\mathfrak{Reg}\left(\mathcal{E}_{\mathcal{\lozenge}_n}\right)$ is a disjoint union of the sets $\mathcal{M}\left(\mathfrak{Reg}\left(\mathcal{E}_{\mathcal{\lozenge}_n}\right)\right)$ and $\mathcal{P}\left(\mathcal{E}_{\mathcal{\lozenge}_n}\right)$.

\vspace{3mm}

\textbf{Proposition 4.3} \textsl{The set $\mathcal{M}\left(\mathfrak{Reg}\left(\mathcal{E}_{\mathcal{\lozenge}_n}\right)\right)$ is a maximal ideal of  semiring $\mathfrak{Reg}\left(\mathcal{E}_{\mathcal{\lozenge}_n}\right)$.}

\emph{Proof.} Obviously,  set  $\mathcal{M}\left(\mathfrak{Reg}\left(\mathcal{E}_{\mathcal{\lozenge}_n}\right)\right)$ is closed under the addition. Let $\alpha \in \mathcal{M}\left(\mathfrak{Reg}\left(\mathcal{E}_{\mathcal{\lozenge}_n}\right)\right)$ and $\beta \in \mathfrak{Reg}\left(\mathcal{E}_{\mathcal{\lozenge}_n}\right)$. Then $(\alpha\cdot \beta)(a_i) = \beta(\alpha(a_i)) = \beta(1) = 1$ for some $a_i$, i.e. $\alpha\cdot \beta \in \mathcal{M}\left(\mathfrak{Reg}\left(\mathcal{E}_{\mathcal{\lozenge}_n}\right)\right)$. Similarly we find $(\beta\cdot\alpha)(a_i) = \alpha(\beta(a_i))$. If $\beta \in \mathcal{M}\left(\mathfrak{Reg}\left(\mathcal{E}_{\mathcal{\lozenge}_n}\right)\right)$, then we choose $a_i$ such that $\beta(a_i) = 1$. If $\beta \in \mathcal{P}\left(\mathcal{E}_{\mathcal{\lozenge}_n}\right)$, we choose $a_i$ such that $\alpha(a_{k_i}) = 1$ where $a_{k_i} = \beta(a_i)$. So, $\beta\cdot \alpha \in \mathcal{M}\left(\mathfrak{Reg}\left(\mathcal{E}_{\mathcal{\lozenge}_n}\right)\right)$. Since every regular endomorphism $\alpha \notin \mathcal{M}\left(\mathfrak{Reg}\left(\mathcal{E}_{\mathcal{\lozenge}_n}\right)\right)$ is an invertible element it follows that $\mathcal{M}\left(\mathfrak{Reg}\left(\mathcal{E}_{\mathcal{\lozenge}_n}\right)\right)$ is a maximal ideal of $\mathfrak{Reg}\left(\mathcal{E}_{\mathcal{\lozenge}_n}\right)$.

\vspace{3mm}

Now we shall extend the maximal ideal $\mathcal{M}\left(\mathfrak{Reg}\left(\mathcal{E}_{\mathcal{\lozenge}_n}\right)\right)$ of $\mathfrak{Reg}\left(\mathcal{E}_{\mathcal{\lozenge}_n}\right)$ to a maximal ideal of  semiring ${\mathcal{E}}_{\mathcal{\lozenge}_n}$. Let us consider all the endomorphisms $\alpha \in \mathcal{E}_{\mathcal{\lozenge}_n}$ such that $\alpha \notin \mathcal{P}\left(\mathcal{E}_{\mathcal{\lozenge}_n}\right)$. Then, if $\alpha \neq \overline{0}$,  we have either $\alpha(1) = a_i$, or $\alpha(a_i) = 1$ for some $i = 1, \ldots, n-2$. We denote the set of all such endomorphisms and $\overline{0}$ by $\mathcal{MAX}\left(\mathcal{E}_{\mathcal{\lozenge}_n}\right)$.

\vspace{3mm}

\textbf{Proposition 4.4} \textsl{For any integer $n > 4$ the set $\mathcal{MAX}\left(\mathcal{E}_{\mathcal{\lozenge}_n}\right)$ is a maximal ideal of semiring ${\mathcal{E}}_{\mathcal{\lozenge}_n}$.}

\emph{Proof.} Let $\alpha, \beta \in \mathcal{MAX}\left(\mathcal{E}_{\mathcal{\lozenge}_n}\right)$, $\alpha \neq \overline{0}$ and $\beta \neq \overline{0}$.

In order to prove that $\mathcal{MAX}\left(\mathcal{E}_{\mathcal{\lozenge}_n}\right)$ is closed under the addition we shall consider four cases.

\emph{Case 1.} Let $\alpha(1) = \beta(1) = a_i$ for some $i = 1, \ldots, n-2$. Then $(\alpha + \beta)(1) = a_i$ and therefore $\alpha + \beta \in \mathcal{MAX}\left(\mathcal{E}_{\mathcal{\lozenge}_n}\right)$.

\emph{Case 2.} Let $\alpha(1) = a_i$ and $\beta(1) = a_j$ where $i, j = 1, \ldots, n-2$ and $i \neq j$. Then $(\alpha + \beta)(1) = 1$. Since $n > 4$, it follows that there are at least three pairwise distinct elements $a_{k_1}$, $a_{k_2}$ and $a_{k_3}$ in  set $\{a_1, \ldots, a_{n-2}\}$. So, even if $\alpha(a_{k_1}) = 0$ and $\beta(a_{k_2}) = 0$, we choose $a_{k_3} \in \{a_1, \ldots, a_{n-2}\}$ such that $\alpha(a_{k_3}) = a_i$ and $\beta(a_{k_3}) = a_j$. Hence, $(\alpha + \beta)(a_{k_3}) = 1$ and therefore $\alpha + \beta \neq i$ and $\alpha + \beta \in \mathcal{MAX}\left(\mathcal{E}_{\mathcal{\lozenge}_n}\right)$.

\emph{Case 3.} Let $\alpha(a_i) = 0$ for some $i = 1, \ldots, n-2$ and $\alpha(1) = 1$. Then $\alpha(a_j) = 1$, where $j = 1, \ldots, n-2$, $j \neq i$, and for any endomorphism $\beta$ it follows $(\alpha + \beta)(1) = 1$ and $(\alpha + \beta)(a_j) = 1$. So, $\alpha + \beta \in \mathcal{MAX}\left(\mathcal{E}_{\mathcal{\lozenge}_n}\right)$.

\emph{Case 4.} Let  $\alpha(a_i) = 1$ for some $i = 1, \ldots, n-2$.  Then it follows $\alpha(1) = 1$. Now for any endomorphism $\beta$ we have $(\alpha + \beta)(1) = 1$ and $(\alpha + \beta)(a_i) = 1$, hence, $\alpha + \beta \in \mathcal{MAX}\left(\mathcal{E}_{\mathcal{\lozenge}_n}\right)$.

In order to prove that $\mathcal{MAX}\left(\mathcal{E}_{\mathcal{\lozenge}_n}\right)\,\cdot\, {\mathcal{E}}_{\mathcal{\lozenge}_n} \subseteq \mathcal{MAX}\left(\mathcal{E}_{\mathcal{\lozenge}_n}\right)$ and ${\mathcal{E}}_{\mathcal{\lozenge}_n}\, \cdot \, \mathcal{MAX}\left(\mathcal{E}_{\mathcal{\lozenge}_n}\right) \subseteq \mathcal{MAX}\left(\mathcal{E}_{\mathcal{\lozenge}_n}\right)$ we shall consider four cases.

\emph{Case 5.} Let $\alpha(1) = a_i$ and $\beta(1) = a_j$ for some $i, j = 1, \ldots, n-2$. Then $(\alpha\cdot\beta)(1) = \beta(\alpha(1)) = \beta(a_i)$. If $\beta(a_i) = 0$, then $\alpha\cdot \beta = \overline{0}$. If $\beta(a_i) = a_j$, it follows $\alpha\cdot \beta \in \mathcal{MAX}\left(\mathcal{E}_{\mathcal{\lozenge}_n}\right)$. By the same reasonings we find $\beta\cdot \alpha \in \mathcal{MAX}\left(\mathcal{E}_{\mathcal{\lozenge}_n}\right)$.

\emph{Case 6.} Let $\alpha(1) = a_i$ and $\beta(1) = 1$ for some $i = 1, \ldots, n-2$. Then $(\alpha\cdot\beta)(1) = \beta(\alpha(1)) = \beta(a_i)$. If $\beta(a_i) = 0$, then $\alpha\cdot \beta = \overline{0}$. If $\beta(a_i) = a_k$, it follows that $\alpha\cdot \beta \in \mathcal{MAX}\left(\mathcal{E}_{\mathcal{\lozenge}_n}\right)$. If $\beta(a_i) = 1$, then $(\alpha\cdot \beta)(1) = 1$. Now if $\alpha(a_i) = 0$, it follows $(\alpha\cdot \beta)(a_i) = 0$ and therefore $\alpha\cdot \beta \in \mathcal{MAX}\left(\mathcal{E}_{\mathcal{\lozenge}_n}\right)$. If $\alpha(a_i) = a_i$, then $(\alpha\cdot \beta)(a_i) = 1$ and also $\alpha\cdot \beta \in \mathcal{MAX}\left(\mathcal{E}_{\mathcal{\lozenge}_n}\right)$.
At last we find  $(\beta\cdot \alpha)(1) = \alpha(\beta(1)) = \alpha(1) = a_i$, so, $ \beta\cdot \alpha \in \mathcal{MAX}\left(\mathcal{E}_{\mathcal{\lozenge}_n}\right)$.

\emph{Case 7.} Let $\alpha(a_i) = 1$ and $\beta(a_j) = 1$ for some $i, j = 1, \ldots, n-2$. Then $(\alpha\cdot \beta)(1) = (\beta\cdot \alpha)(1) = 1$,  $(\alpha\cdot \beta)(a_i) = \beta(\alpha(a_i)) = \beta(1) = 1$ and $(\beta\cdot \alpha)(a_j) + \alpha(\beta(a_j) = \alpha(1) = 1$. Hence, $\alpha\cdot \beta,\, \beta\cdot \alpha \in \mathcal{MAX}\left(\mathcal{E}_{\mathcal{\lozenge}_n}\right)$.

\emph{Case 8.} Let $\alpha(a_i) = 1$ and $\beta(1) = 1$ for some $i = 1, \ldots, n-2$.  Then $(\alpha\cdot\beta)(1) = \beta(\alpha(1)) = \beta(1) = 1$. Since $(\alpha\cdot\beta)(a_i) = \beta(\alpha(a_i)) = \beta(1) = 1$, it follows that $\alpha\cdot \beta \in \mathcal{MAX}\left(\mathcal{E}_{\mathcal{\lozenge}_n}\right)$. Similarly, we have $(\beta\cdot\alpha)(1) = \alpha(\beta(1)) = \alpha(1) = 1$. If we assume that $\beta(a_i) = 1$ for some $a_i$, then $(\beta\cdot\alpha)(a_i) = \alpha(\beta(a_i)) = \alpha(1) = 1$ and follows $\beta\cdot \alpha \in \mathcal{MAX}\left(\mathcal{E}_{\mathcal{\lozenge}_n}\right)$. If we assume that $\beta(a_i) = 0$ for some $a_i$, then $(\beta\cdot\alpha)(a_i) = \alpha(\beta(a_i)) = \alpha(0) = 0$ and it follows $\beta\cdot \alpha \in \mathcal{MAX}\left(\mathcal{E}_{\mathcal{\lozenge}_n}\right)$. Now let $\beta \in \mathcal{P}\left(\mathcal{E}_{\mathcal{\lozenge}_n}\right)$. So, there is $a_j$ such that $\beta(a_j) = a_i$. Then $(\beta\cdot\alpha)(a_j) = \alpha(\beta(a_j)) = \alpha(a_i) = 1$. Therefore in all the possibilities for endomorphism $\beta$ it follows $\beta\cdot \alpha \in \mathcal{MAX}\left(\mathcal{E}_{\mathcal{\lozenge}_n}\right)$.

\vspace{3mm}

From Example 3.2 and the last proposition immediately follows

\vspace{3mm}

\textbf{Corollary 4.5} \textsl{The semiring ${\mathcal{E}}_{\mathcal{\lozenge}_n}$ is a simple semiring only for $n = 4$.}

\vspace{3mm}

Finally, in this section we consider the idempotent endomorphisms of semiring $\mathfrak{Reg}\left(\mathcal{E}_{\mathcal{\lozenge}_n}\right)$. Let us denote the set of such elements by $\mathcal{ID}\left(\mathfrak{Reg}\left(\mathcal{E}_{\mathcal{\lozenge}_n}\right)\right)$. In [ ] we find a result similar to the next proposition.

\vspace{3mm}

\textbf{Proposition 4.6} \textsl{The set $\mathcal{ID}\left(\mathfrak{Reg}\left(\mathcal{E}_{\mathcal{\lozenge}_n}\right)\right)$ is a commutative subsemiring of $\mathfrak{Reg}\left(\mathcal{E}_{\mathcal{\lozenge}_n}\right)$ such that the addition and the multiplication tables coincide.  The subset $I = \{\overline{1}, \psi_{ii}, i = 1, \ldots, n-2\}$ is an ideal of semiring $\mathcal{ID}\left(\mathfrak{Reg}\left(\mathcal{E}_{\mathcal{\lozenge}_n}\right)\right)$.}

\emph{Proof.} Let $\alpha \in \mathcal{ID}\left(\mathfrak{Reg}\left(\mathcal{E}_{\mathcal{\lozenge}_n}\right)\right)$. Since $ 0 \notin \im(\alpha)$, $ \alpha \neq \overline{a_i}$, $\;i = 1, \ldots, n-2$ and $ \alpha(1) = 1$, it follows that $\alpha$ is an identity on some subset $A\subseteq \{a_1, \ldots, a_{n-2}\}$ and $\alpha(x) = \overline{1}$ for any $x \notin A$.\break So, if:

(i) $A = \{a_1, \ldots, a_{n-2}\}$, then $\alpha = i$,

(ii) $A = \varnothing$, then $\alpha = \overline{1}$,

(iii) $A = \{a_i\}$, then $\alpha = \psi_{ii}$, where $i = 1, \ldots, n-2$.

The fact that $\alpha$ is an identity on $A\subseteq \{a_1, \ldots, a_{n-2}\}$ we denote by $\alpha_A$. Now we find that
$$\alpha_A + \alpha_B = \alpha_A\cdot\alpha_B = \alpha_B\cdot\alpha_A = \alpha_{A\cap B}.$$

Hence, $\mathcal{ID}\left(\mathfrak{Reg}\left(\mathcal{E}_{\mathcal{\lozenge}_n}\right)\right)$ is a commutative semiring such that the addition and the multiplication tables coincide. Note that  identity $i$ is  both an additively neutral and  a multiplicatively neutral element of this semiring.

From the last equalities it follows
$$\alpha_A + \psi_{i,i} = \alpha_A\cdot\psi_{i,i} = \psi_{i,i}\cdot\alpha_A = \left\{\begin{array}{cl} \psi_{i,i}, & \; \mbox{if}\; a_i \in A\\ \overline{1}, & \; \mbox{if}\; a_i \notin A \end{array} \right.$$
and $\alpha_A + \overline{1} = \alpha_A\cdot \overline{1} = \overline{1}\cdot \alpha_A = \overline{1}$ for arbitrary $\alpha_A \in \mathcal{ID}\left(\mathfrak{Reg}\left(\mathcal{E}_{\mathcal{\lozenge}_n}\right)\right)$. Hence, $I$ is an ideal of semiring $\mathcal{ID}\left(\mathfrak{Reg}\left(\mathcal{E}_{\mathcal{\lozenge}_n}\right)\right)$.

\vspace{3mm}

\vspace{7mm}

{\bf 5 \hspace{1mm} Endomorphisms with images which are chains}

\vspace{3mm}

Let us consider endomorphisms $\alpha \in {\mathcal{E}}_{\mathcal{\lozenge}_n}$ such that $\im (\alpha) \subseteq \mathcal{C}_i = \{0, a_i, 1\} \subset \mathcal{\lozenge}_n$, where $i = 1, \ldots, n-2$.
For any fixed $i \in \{1, \ldots, n-2\}$ we denote the set of all such endomorphisms $\alpha$ by ${\mathcal{E}}_{\mathcal{\lozenge}_n}(a_i)$.

\vspace{3mm}

\textbf{Proposition 5.1} \textsl{For any $n \geq 4$ and $i \in \{1, \ldots, n-2\}$ the  set ${\mathcal{E}}_{\mathcal{\lozenge}_n}(a_i)$ is a subsemiring with zero of ${\mathcal{E}}_{\mathcal{\lozenge}_n}$.}

\emph{Proof.} It is easy to see that $\im (\alpha_i) \subseteq \mathcal{C}_i$ and $\im (\beta_i) \subseteq \mathcal{C}_i$ imply $\im (\alpha_i + \beta_i) \subseteq \mathcal{C}_i$ and $\im (\alpha_i \cdot \beta_i) \subseteq \mathcal{C}_i$.

\vspace{3mm}

Note that it is easy to prove that for $n \geq 4$ the order of any semiring ${\mathcal{E}}_{\mathcal{\lozenge}_n}(a_i)$, $i = 1, \ldots, n - 2$, is equal to $3(n-1)$.

Let for any fixed $i \in \{1, \ldots, n-2\}$ consider the subset of ${\mathcal{E}}_{\mathcal{\lozenge}_n}(a_i)$ consisting of all the endomorphisms $\alpha$ such that $\alpha(1) \neq 1$ and denote this subset by ${\mathcal{E}}^{[\,0, i\,]}_{\mathcal{\lozenge}_n}$. From Proposition 3.1. a. it follows that $1 \notin \im(\alpha)$ for every $\alpha \in {\mathcal{E}}^{[\,0, i\,]}_{\mathcal{\lozenge}_n}$ and by using Proposition 3.1. b. and notations (3) we find that all the elements of ${\mathcal{E}}^{[\,0, i\,]}_{\mathcal{\lozenge}_n}$ are:
$$\overline{0}= \wr\, 0, \ldots, 0\,\wr,\; \overline{a_i}= \wr\, a_i, \ldots, a_i\,\wr,\;
$$
$$\alpha_{\,01}^{(i)} = \wr\, 0, a_i, \ldots, a_i\, \wr, \; \alpha_{\,02}^{(i)} = \wr\, a_i, 0, a_i, \ldots, a_i\, \wr, \; \cdots, \alpha_{\,0\,n-2}^{(i)} = \wr\, a_i, \ldots, a_i, 0, a_i\, \wr.
$$

It is easy to verify that:
$\overline{0} + \alpha = \alpha$,  $\alpha + \alpha = \alpha$ and $\alpha + \overline{a_i} = \overline{a_i}$ for arbitrary $\alpha \in {\mathcal{E}}^{[\,0, i\,]}_{\mathcal{\lozenge}_n}$ and also $\alpha_{\,0j}^{(i)} + \alpha_{\,0k}^{(i)} = \overline{a_i}$, where $j \neq k$.

 Obviously, $\overline{0}\cdot \alpha = \alpha\cdot \overline{0} = \overline{0}$. By using notations (4) we obtain  $\alpha\cdot \alpha_{0,i} = \overline{0}$ for arbitrary $\alpha \in {\mathcal{E}}^{[\,0, i\,]}_{\mathcal{\lozenge}_n}$. We also find $\alpha\cdot \beta = \alpha$, where $\alpha, \beta \in {\mathcal{E}}^{[\,0, i\,]}_{\mathcal{\lozenge}_n}$ and $\beta \neq \overline{0}$, $\beta \neq \alpha_{0,i}$. Hence, we prove

\vspace{3mm}

\textbf{Proposition 5.2} \textsl{For any $n \geq 4$ and $i \in \{1, \ldots, n-2\}$ the set ${\mathcal{E}}^{[\,0, i\,]}_{\mathcal{\lozenge}_n}$ is a subsemiring of ${\mathcal{E}}_{\mathcal{\lozenge}_n}(a_i)$ with zero element $\overline{0}$.}

\vspace{3mm}

Note that in semiring ${\mathcal{E}}^{[\,0, i\,]}_{\mathcal{\lozenge}_n}$ the endomorphism $\alpha_{0,i}$ is a right zero and all other endomorphisms $\alpha_{0,j}$, where $j \neq i$, are right identities. So, there are not any ideals of order more than two. But for any $n \geq 4$ the semiring ${\mathcal{E}}^{[\,0, i\,]}_{\mathcal{\lozenge}_n}$ is not simple because the two--element set $\{\overline{0}, \alpha_{0i}\}$ is a proper ideal of ${\mathcal{E}}^{[\,0, i\,]}_{\mathcal{\lozenge}_n}$.

\vspace{3mm}

 Now for any fixed $i \in \{1, \ldots, n-2\}$ we consider the subset of ${\mathcal{E}}_{\mathcal{\lozenge}_n}(a_i)$ consisting of all the endomorphisms $\alpha$ such that $0 \notin \im(\alpha)$ and denote this subset by ${\mathcal{E}}^{[\,i, 1\,]}_{\mathcal{\lozenge}_n}$. Let us observe that ${\mathcal{E}}^{[\,i, 1\,]}_{\mathcal{\lozenge}_n}\subseteq {\mathcal{E}}^{[\,a, 1\,]}_{\mathcal{\lozenge}_n}$, see Lemma 3.7.
 By using Proposition 3.1. c. and d. we find that all the elements of ${\mathcal{E}}^{[\,i, 1\,]}_{\mathcal{\lozenge}_n}$ are:
$$ \overline{a_i}= \wr\, a_i, \ldots, a_i\,\wr,\; \overline{1}= \wr\, 1, \ldots, 1\,\wr,\;
$$
$$\psi_{\,1,i} = \wr\, a_i, 1, \ldots, 1\, \wr, \; \psi_{\,2,i} = \wr\, 1, a_i, 1 \ldots, 1\, \wr, \; \cdots, \psi_{\,n-2,i} = \wr\, 1, \ldots, 1, a_i, 1\, \wr. \eqno{(5)}
$$

It is easy to verify that:
$\overline{a_i} + \alpha = \alpha$,  $\alpha + \alpha = \alpha$ and $\alpha + \overline{1} = \overline{1}$ for every $\alpha \in {\mathcal{E}}^{[\,i, 1\,]}_{\mathcal{\lozenge}_n}$ and also $\psi_{j,i} + \psi_{k,i} = \overline{1}$, where $j \neq k$.

We also compute $\alpha\cdot \overline{a_i} = \overline{a_i}$, $\;\alpha\cdot \psi_{i,i} = \alpha$ and $\alpha\cdot \psi_{j,i} = \overline{1}$, where $j \neq i$ for every $\alpha \in {\mathcal{E}}^{[\,i, 1\,]}_{\mathcal{\lozenge}_n}$.  Hence, we prove

\vspace{4mm}

\textbf{Proposition 5.3} \textsl{For any $n \geq 4$ and $i \in \{1, \ldots, n-2\}$ the set ${\mathcal{E}}^{[\,i, 1\,]}_{\mathcal{\lozenge}_n}$ is a subsemiring of ${\mathcal{E}}^{[\,i\,]}_{\mathcal{\lozenge}_n}$ without a zero element.}

\vspace{4mm}

Note that semiring ${\mathcal{E}}^{[\,i, 1\,]}_{\mathcal{\lozenge}_n}$ is not simple. For instance,  set $M = {\mathcal{E}}^{[\,i, 1\,]}_{\mathcal{\lozenge}_n} \backslash \{\psi_{i,i}\}$ is a maximal ideal of this semiring.
\vspace{3mm}

 By analogy with the above reasonings we consider, for any semiring ${\mathcal{E}}_{\mathcal{\lozenge}_n}(a_i)$, the subset of all endomorphisms $\alpha$ such that $\im(\alpha) \subseteq \{0,1\}$. Since these endomorphisms are well defined by (1), we observe that this subset is  semiring ${\mathcal{E}}^{[\,0, 1\,]}_{\mathcal{\lozenge}_n}$. It is straightforward that $\ds {\mathcal{E}}^{[\,0, 1\,]}_{\mathcal{\lozenge}_n} = {\mathcal{E}}_{\mathcal{\lozenge}_n}(a_i)\cap {\mathcal{E}}_{\mathcal{\lozenge}_n}(a_j)$, where $i \neq j$,
 and therefore $\ds {\mathcal{E}}^{[\,0, 1\,]}_{\mathcal{\lozenge}_n} = \bigcap_{i=1}^n {\mathcal{E}}_{\mathcal{\lozenge}_n}(a_i)$.
  So, we have:

 \vspace{4mm}

\textbf{Proposition 5.4} \textsl{For any $n \geq 4$ the  semiring ${\mathcal{E}}^{[\,0, 1\,]}_{\mathcal{\lozenge}_n}$ is a subsemiring  of all semirings ${\mathcal{E}}_{\mathcal{\lozenge}_n}(a_i)$, where $i \in \{1, \ldots, n-2\}$.}

\vspace{3mm}

Note that from the proof of Lemma 3.5 it follows that every nonzero element of  semiring ${\mathcal{E}}^{[\,0, 1\,]}_{\mathcal{\lozenge}_n}$  is a right identity. Hence, ${\mathcal{E}}^{[\,0, 1\,]}_{\mathcal{\lozenge}_n}$  is a simple subsemiring of ${\mathcal{E}}_{\mathcal{\lozenge}_n}(a_i)$.

\vspace{3mm}

By using the  definitions of  semirings from the last three propositions immediately follows

\vspace{3mm}

\textbf{Corollary 5.5} \textsl{For any $n \geq 4$ and $i \in \{1, \ldots, n-2\}$, the semiring  ${\mathcal{E}}_{\mathcal{\lozenge}_n}(a_i)$ is a union of three subsemirings ${\mathcal{E}}^{[\,0, i\,]}_{\mathcal{\lozenge}_n}$, ${\mathcal{E}}^{[\,i, 1\,]}_{\mathcal{\lozenge}_n}$ and ${\mathcal{E}}^{[\,0, 1\,]}_{\mathcal{\lozenge}_n}$.}

 \vspace{3mm}

 Now we prove the main result of the section:

 \vspace{4mm}

\textbf{Theorem 5.6} \textsl{For any $n \geq 4$ and $i \in \{1, \ldots, n-2\}$ the semiring ${\mathcal{E}}_{\mathcal{\lozenge}_n}(a_i)$ is  a simple semiring.}

\vspace{1mm}
\emph{Proof.} We shall use the notations from the proofs of Proposition 4.2, Proposition 4.3, Proposition 4.4 and Lemma 3.5. Let us consider  equalities:
$$  \overline{0}\cdot \psi_{i,i} = \overline{0}, \;\;  \overline{1}\cdot \psi_{i,i} = \overline{1}, \; \:  \varphi_j\cdot \psi_{i,i} = \varphi_j, \;\; \alpha^{(i)}_{0,j}\cdot \psi_{i,i} = \alpha^{(i)}_{0,j}, \; \; \psi_{j,i}\cdot \psi_{i,i} = \psi_{j,i},$$
where $j = 1, \ldots, n-2$. Hence, $\psi_{ii}$ is a right identity of the semiring ${\mathcal{E}}_{\mathcal{\lozenge}_n}(a_i)$.

 Let $J$ be an ideal of  semiring ${\mathcal{E}}_{\mathcal{\lozenge}_n}(a_i)$. Since $\overline{0}$ is a zero element of this semiring, we can suppose that $\overline{0} \in J$. We shall consider the following three cases:

\vspace{2mm}

\emph{Case 1.} Let us assume that $\varphi_1 \in J$. If $i \neq 1$ we find $\varphi_1\cdot \alpha^{(i)}_{0,i} = \alpha^{(i)}_{0,1}\in J$ and $\alpha^{(i)}_{0,i}\cdot \varphi_1 = \varphi_i \in J$. Then $\alpha^{(i)}_{0,1} + \varphi_i = \psi_{i,i} \in J$. Since $\psi_{ii}$ is a right identity, it follows that $J = {\mathcal{E}}_{\mathcal{\lozenge}_n}(a_i)$. If $i = 1$, we find $\varphi_2\cdot \varphi_1 = \varphi_2 \in J$ and $\varphi_2\cdot \alpha^{(1)}_{0,1} = \alpha^{(1)}_{0,2}$. Then $\alpha^{(1)}_{0,2} + \varphi_1 = \psi_{1,1} \in J$ and $J = {\mathcal{E}}_{\mathcal{\lozenge}_n}(a_i)$.

\vspace{2mm}

\emph{Case 2.} Let us assume that $\alpha \in J$, where $\alpha(1) = 1$. Since $\varphi_1\cdot \alpha = \varphi_1 \in J$ we go to Case 1.

\vspace{2mm}

\emph{Case 3.} Let us assume that $\alpha \in J$, where $\alpha(1) = a_i$. Then for $j \neq i$ it follows $\alpha\cdot \varphi_j = \beta \in J$ where $\beta(1) = 1$ and we go to Case 2.

Hence, either $J = \{\overline{0}\}$, or $J = {\mathcal{E}}_{\mathcal{\lozenge}_n}(a_i)$ and this completes the proof.

 \vspace{4mm}

\textbf{Remark 5.7} From the main result (Theorem 1.7) of [9] we know that  endomorphism semiring (with zero) $\mathcal{E}_{\mathcal{C}_n}$ of a finite chain is a simple semiring. About subsemirings, from Theorem 2.1 of [4], it follows that a subsemiring of $\mathcal{E}_{\mathcal{C}_n}$  is simple provided it contains one well known ideal of $\mathcal{E}_{\mathcal{C}_n}$. In [1] we prove that every proper subsemiring of $\mathcal{E}_{\mathcal{C}_n}$ is not simple.

Here we are interested in the  semirings of the endomorphisms of a semilattice $\mathcal{\lozenge}_n$ and prove that

$\bullet$ Semiring ${\mathcal{E}}_{\mathcal{\lozenge}_n}$, when $n > 4$, is not simple.

$\bullet$ Semiring ${\mathcal{E}}_{\mathcal{\lozenge}_n}(a_i)$, for any $i = 1, \ldots, n-2$, is a simple semiring.

$\bullet$ Semiring ${\mathcal{E}}^{[0,1]}_{\mathcal{\lozenge}_n}$  is also a simple semiring.
\vspace{3mm}

 A natural question is: \emph{Is there a semiring $R$ such that ${\mathcal{E}}^{[0,1]}_{\mathcal{\lozenge}_n}\subsetneqq R \subsetneqq {\mathcal{E}}_{\mathcal{\lozenge}_n}(a_i)$ and $R$ is not a simple semiring?}

Let $R$ be a set consisting of  all the elements of semiring ${\mathcal{E}}^{[0,1]}_{\mathcal{\lozenge}_n}$ and also of  all the endomorphisms $\psi_{1,i}, \ldots, \psi_{n-2,i}$, see (5). We observe that $\varphi_k + \psi_{j,i} = \overline{1}$ if $k \neq j$ and $\varphi_j + \psi_{j,i} = \psi_{j,i}$. From the proof of Proposition 5.3. we have $\psi_{j,i} + \psi_{k,i} = \overline{1}$, where $j \neq k$, and $\psi_{j,i} + \overline{1} = \overline{1}$. Then, from the proof of Lemma 3.5 it follows that $R$ is closed under the addition.

Now we find $\varphi_k \cdot \psi_{j,i} = \varphi_k, \; \psi_{j,i} \cdot \varphi_i = \varphi_i, \; \psi_{j,i} \cdot \varphi_k = \overline{1}$ if $ i \neq k$ from the proof of Theorem 3.8. From the equalities $\;\alpha\cdot \psi_{i,i} = \alpha$ and $\alpha\cdot \psi_{j,i} = \overline{1}$, where $j \neq i$ for every $\alpha \in {\mathcal{E}}^{[\,i, 1\,]}_{\mathcal{\lozenge}_n}$, from the proof of Proposition 5.3 and again the proof of Lemma 3.5 it follows that $R$ is closed under the multiplication. Thus we also verify that $R\,{\mathcal{E}}^{[0,1]}_{\mathcal{\lozenge}_n} \subseteq {\mathcal{E}}^{[0,1]}_{\mathcal{\lozenge}_n}$
and ${\mathcal{E}}^{[0,1]}_{\mathcal{\lozenge}_n}\,R \subseteq {\mathcal{E}}^{[0,1]}_{\mathcal{\lozenge}_n}$. So, we prove

 \vspace{4mm}

\textbf{Proposition 5.8} \textsl{The set $R$ is   a  subsemiring of semiring ${\mathcal{E}}_{\mathcal{\lozenge}_n}(a_i)$ which is not simple. The semiring ${\mathcal{E}}^{[0,1]}_{\mathcal{\lozenge}_n}$ is an ideal of $R$.}

 \vspace{3mm}

Using the notations from the begining of  Section 5 we consider the endomorphisms $\alpha \in {\mathcal{E}}_{\mathcal{\lozenge}_n}$ such that $\ds \im(\alpha) \subseteq \bigcup_{\ell=1}^k \mathcal{C}_{i_\ell}$ where $\mathcal{C}_{i_\ell} = \{0, a_{i_\ell}, 1\} \subset \mathcal{\lozenge}_n$, $\ell = 1, \ldots, k$ and $k \leq n-2$. We denote the set of these endomorphisms by ${\mathcal{E}}_{\mathcal{\lozenge}_n}(a_{i_1}, \ldots, a_{i_k})$. After renumbering we can denote this set by ${\mathcal{E}}_{\mathcal{\lozenge}_n}(a_{1}, \ldots, a_{k})$, where $k \leq n-2$. By similar reasonings to those used in the proof  of Proposition 5.1 we prove

 \vspace{4mm}

\textbf{Proposition 5.9} \textsl{For any $1 \leq k \leq n-2$ the set ${\mathcal{E}}_{\mathcal{\lozenge}_n}(a_{1}, \ldots, a_{k})$ is a  subsemiring with zero of semiring ${\mathcal{E}}_{\mathcal{\lozenge}_n}$.}

\vspace{8mm}

{\bf 6 \hspace{1mm} Idempotent elements of  the endomorphism semiring}

\vspace{3mm}

The set of all the idempotents of  semiring ${\mathcal{E}}_{\mathcal{\lozenge}_n}$ is not closed under the multiplication. For instance,  endomorphisms $\wr\, 0, 1, \ldots, 1\,\wr$ and $\wr a_1, 0, a_1, \ldots, a_1\,\wr$ are idempotents, but their product $\wr\, 0, 1, \ldots, 1\,\wr \cdot\wr a_1, 0, a_1, \ldots, a_1\,\wr = \wr 0, a_1, \ldots, a_1\,\wr$ is not an idempotent. So, we shall consider only a part of the idempotent elements of ${\mathcal{E}}_{\mathcal{\lozenge}_n}$.

\vspace{3mm}

Idempotent elements $\alpha \in {\mathcal{E}}_{\mathcal{\lozenge}_n}$ such that  $\alpha(1) = 1$ are called \emph{stable idempotents}. The set of idempotents which are not stable is not closed under the addition and multiplication. For instance, $\wr a_1, 0, a_1, \ldots, a_1\,\wr$ and
$\wr a_2, a_2, 0, a_2, \ldots, a_2\,\wr$ are idempotents which are not stable but
$$\wr a_1, 0, a_1, \ldots, a_1\,\wr + \wr a_2, a_2, 0, a_2, \ldots, a_2\,\wr = \wr 1, a_2, a_1, 1, \ldots, 1\,\wr\; \mbox{and}$$
$$\wr a_1, 0, a_1, \ldots, a_1\,\wr \cdot \wr a_2, a_2, 0, a_2, \ldots, a_2\,\wr = \wr a_2, 0, a_2, \ldots, a_2\,\wr$$
are not idempotent endomorphisms.

The set of stable idempotent elements of semiring ${\mathcal{E}}_{\mathcal{\lozenge}_n}$ we denote by $\mathcal{SI}\left(\mathcal{E}_{\mathcal{\lozenge}_n}\right)$. Let $\alpha \in \mathcal{SI}\left(\mathcal{E}_{\mathcal{\lozenge}_n}\right)$ and $\alpha(a_i) = 0$. Since $\alpha(1) = 1$ from Proposition 3.1, b., we obtain that $\alpha(a_j) = 1$ for any $j \neq i$. So, $\alpha = \varphi_i$. Thus, it follows that $\left({\mathcal{E}}^{[0,1]}_{\mathcal{\lozenge}_n}\right)^* \subset \mathcal{SI}\left(\mathcal{E}_{\mathcal{\lozenge}_n}\right)$.

\vspace{3mm}

\textbf{Lemma 6.1} \textsl{The endomorphism $\alpha \in \mathcal{SI}\left(\mathcal{E}_{\mathcal{\lozenge}_n}\right) \backslash \left({\mathcal{E}}^{[0,1]}_{\mathcal{\lozenge}_n}\right)^*$ if and only if either $\alpha(a_i) = a_i$, or $\alpha(a_i) = 1$ for $i = 1, \ldots, n-2$.}

\emph{Proof.} Let $\alpha \in {\mathcal{E}}_{\mathcal{\lozenge}_n}$ and either $\alpha(a_i) = a_i$, or $\alpha(a_i) = 1$ for $i = 1, \ldots, n-2$. Then either $\alpha^2(a_i) = \alpha(a_i) = a_i$, or $\alpha^2(a_i) = \alpha(1) = 1$ that is $\alpha^2 = \alpha$.

Conversely, let $\alpha \in \mathcal{SI}\left(\mathcal{E}_{\mathcal{\lozenge}_n}\right) \backslash \left({\mathcal{E}}^{[0,1]}_{\mathcal{\lozenge}_n}\right)^*$. We assume that $\alpha(a_i) \neq 1$. Then $\alpha(a_i) = a_k$ and $\alpha^2(a_i) = \alpha(a_k)$. But $\alpha^2 = \alpha$. Hence, $\alpha^2(a_i) = \alpha(a_i)$. Thus, it follows $\alpha(a_k) = \alpha^2(a_i) = \alpha(a_i) = a_k$. From Proposition 3.1, d., we obtain that $\alpha$ is a permutation, that is $\alpha(a_i) = \alpha(a_k)$ means that $a_i = a_k$. So either $\alpha(a_i) = 1$, or $\alpha(a_i) = a_i$.

\vspace{3mm}

Immediately from the last lemma and Proposition 4.6 it follows that
$$\mathcal{SI}\left(\mathcal{E}_{\lozenge_n}\right) =
\mathcal{ID}\left(\mathfrak{Reg}\left(\mathcal{E}_{\mathcal{\lozenge}_n}\right)\right) \cup \left({\mathcal{E}}^{[0,1]}_{\mathcal{\lozenge}_n}\right)^*.$$

\textbf{Proposition 6.2} \textsl{The set $\mathcal{SI}\left(\mathcal{E}_{\lozenge_n}\right)$ is a subsemiring of ${\mathcal{E}}_{\mathcal{\lozenge}_n}$. Set $\left({\mathcal{E}}^{[0,1]}_{\mathcal{\lozenge}_n}\right)^*$ is an ideal of $\mathcal{SI}\left(\mathcal{E}_{\mathcal{\lozenge}_n}\right)$. The set $\mathcal{SI}\left(\mathcal{E}_{\lozenge_n}\right)\backslash \{i\}$ is a maximal ideal of $\mathcal{SI}\left(\mathcal{E}_{\lozenge_n}\right)$.}

\emph{Proof.} For the first part of the theorem it is enough to show that if $\alpha \in \mathcal{ID}\left(\mathfrak{Reg}\left(\mathcal{E}_{\mathcal{\lozenge}_n}\right)\right)$ and $\beta \in
\left({\mathcal{E}}^{[0,1]}_{\mathcal{\lozenge}_n}\right)^*$, then $\alpha + \beta \in \mathcal{SI}\left(\mathcal{E}_{\lozenge_n}\right)$, $\alpha\cdot \beta \in \left({\mathcal{E}}^{[0,1]}_{\mathcal{\lozenge}_n}\right)^*$ and $\beta\cdot \alpha \in \left({\mathcal{E}}^{[0,1]}_{\mathcal{\lozenge}_n}\right)^*$. By the notations from the proof of Propsition 4.6 we choose $\alpha = \alpha_A$, where $A\subseteq \{a_1, \ldots, a_{n-2}\}$. Now it follows
$\ds \varphi_i + \alpha_A  = \left\{\begin{array}{cl} \psi_{i,i}, & \; \mbox{if}\; a_i \in A\\ \overline{1}, & \; \mbox{if}\; a_i \notin A \end{array} \right.$,
where $i = 1, \ldots, n-2$. Therefore, $\mathcal{SI}\left(\mathcal{E}_{\lozenge_n}\right)$ is closed under the addition. From $\ds  \alpha_A\cdot \varphi_i  = \left\{\begin{array}{cl} \varphi_i, & \; \mbox{if}\; a_i \in A\\ \overline{1}, & \; \mbox{if}\; a_i \notin A \end{array} \right.$ and $\varphi_i\cdot \alpha_A = \varphi_i$ it follows that $\mathcal{SI}\left(\mathcal{E}_{\lozenge_n}\right)$ is a semiring and $\left({\mathcal{E}}^{[0,1]}_{\mathcal{\lozenge}_n}\right)^*$ is an ideal of this semiring. From the above reasonings it follows that there are not any endomorphisms $\alpha, \beta \in \mathcal{SI}\left(\mathcal{E}_{\lozenge_n}\right)$, $\alpha \neq i$ and $\beta \neq i$, such that $\alpha + \beta = i$ or $\alpha\cdot \beta = i$. Hence, the set $\mathcal{SI}\left(\mathcal{E}_{\lozenge_n}\right)\backslash \{i\}$ is a maximal ideal of semiring $\mathcal{SI}\left(\mathcal{E}_{\lozenge_n}\right)$.

\vspace{3mm}

The semiring $\mathcal{SI}\left(\mathcal{E}_{\lozenge_n}\right)$ is an extension of the commutative semiring $\mathcal{ID}\left(\mathfrak{Reg}\left(\mathcal{E}_{\mathcal{\lozenge}_n}\right)\right)$ in the set of all the idempotent elements of  semiring ${\mathcal{E}}_{\mathcal{\lozenge}_n}$ and this extension is a non-commutative semiring without zero.
Now we shall extend the semiring of regular idempotents $\mathcal{ID}\left(\mathfrak{Reg}\left(\mathcal{E}_{\mathcal{\lozenge}_n}\right)\right)$ to another non-commutative semiring with zero whose elements are idempotents by the following reasonings. After Proposition 3.3 we note that the elements of the right ideal $\mathcal{AC}\left({\mathcal{E}}_{\mathcal{\lozenge}_n}\right)$  of $\mathcal{E}_{\mathcal{\lozenge}_n}$ are idempotents. Now we consider  set $\widehat{\mathcal{ID}}\left(\mathfrak{Reg}\left(\mathcal{E}_{\mathcal{\lozenge}_n}\right)\right) = \mathcal{ID}\left(\mathfrak{Reg}\left(\mathcal{E}_{\mathcal{\lozenge}_n}\right)\right)\cup \mathcal{AC}\left({\mathcal{E}}_{\mathcal{\lozenge}_n}\right)$.

\vspace{3mm}

\textbf{Proposition 6.3} \textsl{The set $\widehat{\mathcal{ID}}\left(\mathfrak{Reg}\left(\mathcal{E}_{\mathcal{\lozenge}_n}\right)\right)$ is a  subsemiring with zero of $\mathcal{E}_{\mathcal{\lozenge}_n}$. The right ideal $\mathcal{AC}\left({\mathcal{E}}_{\mathcal{\lozenge}_n}\right)$ is an ideal of this semiring.}

\emph{Proof.} It is enough to show that if $\alpha \in \mathcal{ID}\left(\mathfrak{Reg}\left(\mathcal{E}_{\mathcal{\lozenge}_n}\right)\right)$ and $\beta \in
\mathcal{AC}\left({\mathcal{E}}_{\mathcal{\lozenge}_n}\right)$, then $\alpha + \beta \in \widehat{\mathcal{ID}}\left(\mathfrak{Reg}\left(\mathcal{E}_{\mathcal{\lozenge}_n}\right)\right)$, $\alpha\cdot \beta \in \mathcal{AC}\left({\mathcal{E}}_{\mathcal{\lozenge}_n}\right)$ and $\beta\cdot \alpha \in \mathcal{AC}\left({\mathcal{E}}_{\mathcal{\lozenge}_n}\right)$. By notations from the proof of Propsition 4.6 we choose $\alpha = \alpha_A$, where $A\subseteq \{a_1, \ldots, a_{n-2}\}$. Now it follows
$\ds \overline{a_i} + \alpha_A  = \left\{\begin{array}{cl} \psi_{i,i}, & \; \mbox{if}\; a_i \in A\\ \overline{1}, & \; \mbox{if}\; a_i \notin A \end{array} \right.$,
where $i = 1, \ldots, n-2$. So, $\widehat{\mathcal{ID}}\left(\mathfrak{Reg}\left(\mathcal{E}_{\mathcal{\lozenge}_n}\right)\right)$ is closed under the addition. From $\ds \overline{a_i}\cdot \alpha_A  = \left\{\begin{array}{cl} \overline{a_i}, & \; \mbox{if}\; a_i \in A\\ \overline{1}, & \; \mbox{if}\; a_i \notin A \end{array} \right.$ and $\alpha_A\cdot \overline{a_i} = \overline{a_i}$ follows it that $\widehat{\mathcal{ID}}\left(\mathfrak{Reg}\left(\mathcal{E}_{\mathcal{\lozenge}_n}\right)\right)$ is a semiring and $\mathcal{AC}\left({\mathcal{E}}_{\mathcal{\lozenge}_n}\right)$ is an ideal of this semiring. Obviously, $\overline{0}$ is a zero element of $\widehat{\mathcal{ID}}\left(\mathfrak{Reg}\left(\mathcal{E}_{\mathcal{\lozenge}_n}\right)\right)$.

\vspace{5mm}

In the rest of this section we shall examine the subsemirings of $\mathcal{SI}\left(\mathcal{E}_{\lozenge_n}\right)$.

\vspace{3mm}

\textbf{Example 6.4} Let us consider the three element set $\{\varphi_i, \psi_{i,i}, \overline{1}\}$ where
 $$\varphi_i = \wr\, 1, \ldots, 1, 0, 1, \ldots, 1\,\wr, \psi_{i,i} = \wr\, 1, \ldots, 1, a_i, 1, \ldots, 1\,\wr \; \mbox{and} \; \overline{1} = \wr\, 1, \ldots,  1\,\wr.$$

It is easy to see that this set is the semiring $\mathcal{SI}\left(\mathcal{E}_{\lozenge_n}\right)\cap {\mathcal{E}}_{\mathcal{\lozenge}_n}(a_i)$ which we denote by $\mathcal{SI}(a_i)$. The addition and multiplication tables of semiring $\mathcal{SI}(a_i)$ are:
$$\begin{array}{l|lll}
+ & \varphi_i & \psi_{i,i} & \overline{1}\\ \hline
\varphi_i & \varphi_i & \psi_{i,i} & \overline{1}\\
\psi_{i,i} & \psi_{i,i} & \psi_{i,i} & \overline{1}\\
\overline{1} & \overline{1} & \overline{1} & \overline{1}
\end{array}\quad,\qquad \qquad \begin{array}{l|lll}
\cdot & \varphi_i & \psi_{i,i} & \overline{1}\\ \hline
\varphi_i & \varphi_i & \varphi_i & \varphi_i\\
\psi_{i,i} & \varphi_i & \psi_{i,i} & \overline{1}\\
\overline{1} & \overline{1} & \overline{1} & \overline{1}
\end{array}\;.$$

There are not zero and $\infty$ in this semiring, the endomorphism $\psi_{i,i}$ is an identity and the set $\{\varphi_i,\overline{1}\}$ is an ideal of $\mathcal{SI}(a_i)$.

\vspace{3mm}

Surprisingly, this semiring is an example in [7] of a simple semiring of order 3.

\vspace{3mm}

Clearly, there are $n-2$ three element semirings $\mathcal{SI}(a_i)$, where $i = 1, \ldots, n-2$, and they are isomorphic.
Similarly, we construct  semiring $\mathcal{SI}(a_{i_1}, \ldots, a_{i_k}) = \mathcal{SI}\left(\mathcal{E}_{\lozenge_n}\right)\cap  {\mathcal{E}}_{\mathcal{\lozenge}_n}(a_{i_1}, \ldots, a_{i_k})$.

\vspace{3mm}

\textbf{Proposition 6.5} \textsl{In the semiring $\mathcal{SI}(a_{i_1}, \ldots, a_{i_k})$ there are not a zero element and an element $\infty$ but there is an identity $i(a_{i_1}, \ldots, a_{i_k})$. The set $J = \mathcal{SI}(a_{i_1}, \ldots, a_{i_k}) \backslash \{i(a_{i_1}, \ldots, a_{i_k})\}$ is a maximal ideal of $\mathcal{SI}(a_{i_1}, \ldots, a_{i_k})$.}

\vspace{1mm}

\emph{Proof.} Let $\alpha(a_{i_\ell}) = a_{i_\ell}$ for all $\ell = 1, \ldots, k$ and $\alpha(x) = 1$ for $x \notin \{a_{i_1}, \ldots, a_{i_k}\}$. We denote this endomorphism $\alpha$ by $i(a_{i_1}, \ldots, a_{i_k})$. It  follows easily that $i(a_{i_1}, \ldots, a_{i_k}) \in \mathcal{SI}(a_{i_1}, \ldots, a_{i_k})$. For arbitrary $\varphi_i$, $i = 1, \ldots, n-2$ we have $\varphi_i\cdot i(a_{i_1}, \ldots, a_{i_k}) = \varphi_i$. Let $\alpha \in \mathcal{SI}(a_{i_1}, \ldots, a_{i_k})$ and $\alpha \notin \left({\mathcal{E}}^{[0,1]}_{\mathcal{\lozenge}_n}\right)^*$. Then either $\alpha(a_{i_\ell}) = a_{i_\ell}$, or $\alpha(a_{i_\ell}) = \overline{1}$ where $1 \leq \ell \leq k$.

Since $(\alpha\cdot i(a_{i_1}, \ldots, a_{i_k}))(a_{i_\ell}) = i(a_{i_1}, \ldots, a_{i_k})(\alpha(a_{i_\ell})) = \alpha(a_{i_\ell})$, it follows that $\alpha\cdot i(a_{i_1}, \ldots, a_{i_k}) = \alpha$. So, $i(a_{i_1}, \ldots, a_{i_k})$ is a right identity of semiring $\mathcal{SI}(a_{i_1}, \ldots, a_{i_k})$. Since for any $i = 1, \ldots, n-2$ the endomorphisms $\varphi_i$ are right zeroes, it follows that there are not a zero element and an infinity in this semiring.

Now we compute $\varphi_{i_\ell} + i(a_{i_1}, \ldots, a_{i_k}) = \psi_{\,i_\ell,i_\ell}$ and $\varphi_{i_\ell} + i(a_{i_1}, \ldots, a_{i_k}) = \overline{1}$, where $j \neq i_\ell$, for $\ell = 1, \ldots, k$ and also $\alpha + i(a_{i_1}, \ldots, a_{i_k}) = \alpha$ for all the other elements $\alpha \in \mathcal{SI}(a_{i_1}, \ldots, a_{i_k})$. Let $\alpha, \beta \in \mathcal{SI}(a_{i_1}, \ldots, a_{i_k})$ and $\alpha \neq i(a_{i_1}, \ldots, a_{i_k})$, $\beta \neq i(a_{i_1}, \ldots, a_{i_k})$. Then for some element $a_{i_\ell}$ we have $\alpha(a_{i_\ell}) = 1$. This implies $(\alpha\cdot \beta)(a_{i_\ell}) = 1$, so, $\alpha\cdot \beta \neq i(a_{i_1}, \ldots, a_{i_k})$. Hence, $J = \mathcal{SI}(a_{i_1}, \ldots, a_{i_k}) \backslash \{i(a_{i_1}, \ldots, a_{i_k})\}$ is a maximal ideal of semiring $\mathcal{SI}(a_{i_1}, \ldots, a_{i_k})$.

\vspace{8mm}

{\bf 7 \hspace{1mm} Many examples of simple semirings}

\vspace{3mm}

First, we shall try to find some  ``small'' simple subsemirings of ${\mathcal{E}}_{\mathcal{\lozenge}_n}(a,b)$.

\vspace{3mm}

Let us denote $S_{n-1} =  \left({\mathcal{E}}^{[0,1]}_{\mathcal{\lozenge}_n}\right)^* = {\mathcal{E}}^{[0,1]}_{\mathcal{\lozenge}_n}\backslash \{\overline{0}\}$ which is, obviously, a subsemiring of ${\mathcal{E}}_{\mathcal{\lozenge}_n}$ without a zero element. Using the notations from Example 3.2 it follows that the addition and multiplication tables of semiring $S_3$ are:
$$\begin{array}{c|ccc}
+ & \wr\,011\,\wr & \wr\,101\,\wr  & \overline{1}\\ \hline
\wr\,011\,\wr & \wr\,011\,\wr & \overline{1} & \overline{1}\\
\wr\,101\,\wr & \overline{1} & \wr\,101\,\wr & \overline{1}\\
\overline{1} & \overline{1} & \overline{1} & \overline{1}
\end{array}\quad,\qquad \qquad \begin{array}{c|ccc}
\cdot & \wr\,011\,\wr & \wr\,101\,\wr & \overline{1}\\ \hline
\wr\,011\,\wr & \wr\,011\,\wr & \wr\,011\,\wr & \wr\,011\,\wr\\
\wr\,101\,\wr & \wr\,101\,\wr & \wr\,101\,\wr & \wr\,101\,\wr\\
\overline{1} & \overline{1} & \overline{1} & \overline{1}
\end{array}\;.$$

Since every element of $S_3$ is a right identity, it follows that $S_3$ is a simple semiring without a zero element. Let us consider the endomorphisms $\alpha \in S_3$ such that $\alpha(a) = 1$ (using the notations from Example 3.2). Since $\alpha(1) = 1$, it follows that  two endomorphisms with this property formed a semiring which we denote by $S_2$. The addition and multiplication tables are:
$$\begin{array}{c|cc}
+ &  \wr\,101\,\wr  & \overline{1}\\ \hline
\wr\,101\,\wr & \wr\,101\,\wr &  \overline{1}\\
\overline{1} & \overline{1} & \overline{1}
\end{array}\quad,\qquad \qquad \begin{array}{c|cc}
\cdot &  \wr\,101\,\wr & \overline{1}\\ \hline
\wr\,101\,\wr & \wr\,101\,\wr & \wr\,101\,\wr\\
\overline{1} &  \overline{1} & \overline{1}
\end{array}\;.$$

The same arguments prove that the semiring $S_2$ is a simple semiring without zero.

\vspace{3mm}

Similarly, from the equalities $\varphi_i \cdot \varphi_j = \varphi_i$ and $\overline{1}\cdot\varphi_j = \overline{1}$ for any $i, j \in \{1, \ldots, n-2\}$
from the proof of Lemma 3.5, it follows that  semiring $S_n$ is a simple semiring without zero for any $n \geq 4$. Hence, we can construct a chain of simple semirings without zero
$$S_2 \subset S_3 \subset \cdots \subset S_n.$$
Of course, here $S_2$ and $S_3$ are not the semirings from the examples above, but they are semirings isomorphic to them, respectively.

\vspace{3mm}

Note that the similar construction is realized if we consider the semiring ${\mathcal{E}}^{[0,i]}_{\mathcal{\lozenge}_n}\backslash \{\overline{0},\alpha^{(i)}_{0,i}\}$ which is isomorphic to $S_{n-2}$.

The semirings ${\mathcal{E}}^{[i,1]}_{\mathcal{\lozenge}_n}$ gives other examples of simple semirings but they are trivial. Indeed, the semiring
${\mathcal{E}}^{[i,1]}_{\mathcal{\lozenge}_n}\backslash \{\overline{a_i}, \psi_{i,i}\}$ is a simple and not isomorphic to $S_{n-2}$, but the multiplication is trivial since $\alpha\cdot\beta = \overline{1}$ for all $\alpha$ and $\beta$ of this semiring.

 \vspace{3mm}

Now we shall show that there are ``big'' subsemirings of ${\mathcal{E}}_{\mathcal{\lozenge}_n}$ which are simple.

 \vspace{3mm}

As a consequence of Proposition 5.9 it follows that ${\mathcal{E}}_{\mathcal{\lozenge}_n}(a,b)$, where $a, b \in \{a_1, \ldots , a_{n-2}\}$ is a semiring.

 \vspace{3mm}

\textbf{Theorem 7.1} \textsl{Let $a, b \in \{a_1, \ldots, a_{n-2}\}$. Then for any $n \geq 5$ the semiring ${\mathcal{E}}_{\mathcal{\lozenge}_n}(a,b)$ is a simple  subsemiring of  ${\mathcal{E}}_{\mathcal{\lozenge}_n}$.}

\vspace{1mm}

\emph{Proof.} First, in a similar way, as in the proof of Theorem 5.6, we consider  endomorphism $i(a,b) = \wr\,a, b, 1, \ldots, 1\,\wr$. For arbitrary $\alpha \in {\mathcal{E}}_{\mathcal{\lozenge}_n}(a,b)$ we have $\alpha(x) \in \{0, a, b, 1\}$ for any $x \in \mathcal{\lozenge}_n$. Then $(\alpha\cdot i(a,b))(x) = i(a,b)(\alpha(x)) = \alpha(x)$. Hence, $\alpha\cdot i(a,b) = \alpha$, that is $i(a,b)$ is a right identity of semiring ${\mathcal{E}}_{\mathcal{\lozenge}_n}(a,b)$.

 Let $J$ be an ideal of  semiring ${\mathcal{E}}_{\mathcal{\lozenge}_n}(a,b)$. Since $\overline{0}$ is a zero element of ${\mathcal{E}}_{\mathcal{\lozenge}_n}(a,b)$, we can suppose that $\overline{0} \in J$. We shall consider the following four cases:

\vspace{1mm}

\emph{Case 1.} Let us assume that $\varphi_1 = \wr\, 0, 1, \ldots, 1\, \wr \in J$. For $\alpha^{(2)}_{0,2} = \wr \,b, 0, b, \ldots, b\, \wr$ we calculate $\varphi_1\cdot \alpha^{(2)}_{0,2} = \alpha^{(2)}_{0,1} = (\wr\, 0, b, \ldots, b\, \wr \in J$ and $\alpha^{(2)}_{0,2}\cdot \varphi_1 = \varphi_2 = \wr\, 1, 0, 1, \ldots, 1\, \wr \in J$. Now for $\alpha^{(1)}_{0,1} = \wr\, 0, a, \ldots, a\, \wr$ we find $\varphi_2\cdot \alpha^{(1)}_{0,1} = \alpha^{(1)}_{0,2} = \wr \,a, 0, a, \ldots, a\, \wr \in J$. Hence, it follows $\alpha^{(1)}_{0,2} + \alpha^{(2)}_{0,1} = i(a,b) \in J$ which means that $J = {\mathcal{E}}_{\mathcal{\lozenge}_n}(a,b)$.

\vspace{1mm}

\emph{Case 2.} Let us assume that $\alpha \in J$, where $\alpha(1) = 1$. Since $\varphi_1\cdot \alpha = \varphi_1 \in J$ we go to Case 1.

\vspace{1mm}

\emph{Case 3.} Let us assume that $\alpha \in J$, where $\alpha(1) = a$. Then  $\alpha\cdot \varphi_2 = \beta \in J$, where $\beta(1) = 1$, and we go to Case 2.

\vspace{1mm}

\emph{Case 4.} Let us assume that $\alpha \in J$, where $\alpha(1) = b$. Then  $\alpha\cdot \varphi_1 = \beta \in J$, where $\beta(1) = 1$, and we go to Case 2.

\vspace{1mm}

Hence, either $J = \{\overline{0}\}$, or $J = {\mathcal{E}}_{\mathcal{\lozenge}_n}(a,b)$ and this completes the proof.

\vspace{3mm}

Note that ${\mathcal{E}}_{\mathcal{\lozenge}_n}(a)$ is a left ideal of ${\mathcal{E}}_{\mathcal{\lozenge}_n}(a,b)$.

 \vspace{3mm}

\textbf{Remark 7.2} The semiring ${\mathcal{E}}_{\mathcal{\lozenge}_n}(a,b)$ has a subsemiring, not included in ${\mathcal{E}}_{\mathcal{\lozenge}_n}(a)$ or in ${\mathcal{E}}_{\mathcal{\lozenge}_n}(b)$, which is not simple. For instance, let $n = 5$ and let $S$ be a set of endomorphisms:
$$\overline{1}, \; \psi_{1,a} = \wr\, a, 1, 1, 1\,\wr, \; \psi_{2,a} = \wr\, 1, a, 1, 1\,\wr, \; \psi_{3,a} = \wr\, 1, 1, a, 1\,\wr, \; \psi_{1,b} = \wr\, b, 1, 1, 1\,\wr, \; \psi_{2,b} = \wr\, 1, b, 1, 1\,\wr,$$ $$ \psi_{3,b} = \wr\, 1, 1, b, 1\,\wr, \; \wr\, a, b, 1, 1\,\wr, \; \wr\, b, a, 1, 1\,\wr, \; \wr\, a, 1, b, 1\,\wr, \;\wr\, b, 1, a, 1\,\wr, \; \wr\, 1, a, b, 1\,\wr, \; \wr\, 1, b, a, 1\,\wr.$$

It is easy to establish that $S$ is a subsemiring of ${\mathcal{E}}_{\mathcal{\lozenge}_5}(a,b)$ and $I = \{\overline{1},\psi_{3,a},\psi_{3,b}\}$ is an ideal of $S$. Note also that the idempotent elements of $S$ form commutative semiring such that the addition and the multiplication tables coincide.

\vspace{3mm}

It is intersting to know is there a simple subsemiring of ${\mathcal{E}}_{\mathcal{\lozenge}_n}(a,b)$ which is not included in ${\mathcal{E}}_{\mathcal{\lozenge}_n}(a)$ or in ${\mathcal{E}}_{\mathcal{\lozenge}_n}(b)$. Now we will answer to this question.

 \vspace{3mm}

\textbf{Proposition 7.3} \textsl{For any $n > 4$ there is a subsemiring of ${\mathcal{E}}_{\mathcal{\lozenge}_n}(a,b)$ isomorphic to ${\mathcal{E}}_{\mathcal{\lozenge}_4}$.}

\emph{Proof.} Let, like in Example 3.2, us denote the elements of the lattice ${\mathcal{\lozenge}_4}$ by $0, a, b$ and $1$. Without loss of generality we may suppose that elements of $\mathcal{\lozenge}_n$ are $0, a, b, a_3, \ldots, a_{n-2}$ and $1$. Now we construct a map $\Phi : {\mathcal{E}}_{\mathcal{\lozenge}_4} \rightarrow {\mathcal{E}}_{\mathcal{\lozenge}_n}(a,b)$ such that for any $\alpha \in {\mathcal{E}}_{\mathcal{\lozenge}_4}$, if we denote
$$ \alpha = \wr\, \alpha(a), \alpha(b), \alpha(1)\,\wr, \; \; \mbox{then}\; \; \Phi(\alpha) = \wr\, \alpha(a), \alpha(b), \alpha(1), \ldots, \alpha(1)\,\wr.$$

Since $\alpha(a)$, $\alpha(b)$ and $\alpha(1)$ are elements of ${\mathcal{\lozenge}_4}$, it follows that $\Phi(\alpha) \in {\mathcal{E}}_{\mathcal{\lozenge}_n}(a,b)$. It is easy to obtain that $\Phi(\alpha + \beta) = \Phi(\alpha) + \Phi(\beta)$ and $\Phi(\alpha\cdot \beta) = \Phi(\alpha)\cdot \Phi(\beta)$. Hence, $\phi\left({\mathcal{\lozenge}_4}\right)$ is a subsemiring of ${\mathcal{E}}_{\mathcal{\lozenge}_n}(a,b)$ isomorphic to ${\mathcal{E}}_{\mathcal{\lozenge}_4}$.

 \vspace{3mm}

 A direct consequence of the last proposition is:

 \vspace{3mm}

\textbf{Corollary 7.4} \textsl{For any $n > 4$ semiring ${\mathcal{E}}_{\mathcal{\lozenge}_n}(a,b)$ has a simple subsemiring of order 16.}

\vspace{3mm}

Now we return to semiring ${\mathcal{E}}_{\mathcal{\lozenge}_n}(a_{1}, \ldots, a_{k})$ where $k$ is more than 2. Let us consider  sets $A = \{a_{1}, \ldots, a_{k}\}$ and $A_j = A \backslash \{a_{j}\}$ for $j = 1, \ldots, k$. We denote by ${\mathcal{E}}_{\mathcal{\lozenge}_n}(B)$  semiring ${\mathcal{E}}_{\mathcal{\lozenge}_n}(b_1, \ldots, b_s)$, where $B = \{b_1, \ldots, b_s\}$ is a subset of $A$. So,  semirings ${\mathcal{E}}_{\mathcal{\lozenge}_n}(A)$ and ${\mathcal{E}}_{\mathcal{\lozenge}_n}(A_j)$, $j = 1, \ldots, k$ are well defined.

 \vspace{3mm}

By similar reasonings, as in the proof of Proposition 4.4, we prove the next theorem.

 \vspace{3mm}

\textbf{Theorem 7.5} \textsl{For any integer $k$, where $2 < k \leq n -2$, the set $\ds I = \bigcup_{j=1}^k {\mathcal{E}}_{\mathcal{\lozenge}_n}(A_j)$ is a maximal ideal of semiring ${\mathcal{E}}_{\mathcal{\lozenge}_n}(A)$.}

\vspace{2mm}

\emph{Proof.} All the elements $\alpha \in {\mathcal{E}}_{\mathcal{\lozenge}_n}(A)$ have the property that either $\alpha(1) = a_{\ell}$, where $\ell = 1, \ldots, k$, or $\alpha(1) = 1$ and there are at least $n - k -2$ elements $a_i$ such that $\alpha(a_i) = 1$. Note that  set $\{a_1, \ldots, a_{n-2}\}\backslash A$ is included in the set of all these elements $a_i$. So, the representation of $\alpha$ as an ordered $n-1$ -- tuple consists of more than $n-k-1$ coordinates equal to 1. Analogously, the elements of  set $I$ have the similar property, but these endomorphisms transform at least $n - k$ elements from $\mathcal{\lozenge}_n$ to $1$.

In order to prove that  set $I$ is closed under the addition we shall consider four cases.  Let $\alpha, \beta \in I$, $\alpha \neq \overline{0}$ and $\beta \neq \overline{0}$.
\vspace{1mm}

\emph{Case 1.} Let $\alpha(1) = \beta(1) = a_i$ for some $i = 1, \ldots, k$. Then $(\alpha + \beta)(1) = a_i$ and therefore $\alpha + \beta \in I$.

\vspace{1mm}

\emph{Case 2.} Let $\alpha(1) = a_i$ and $\beta(1) = a_j$ where $i, j = 1, \ldots, k$ and $i \neq j$. Then $(\alpha + \beta)(1) = 1$. Here there are two possibilities. The first one is if $\alpha = \alpha^{(i)}_{0,p}$ and $\beta =  \alpha^{(j)}_{0,p}$, where $p = 1, \ldots, n-2$. Then $\alpha + \beta$ transforms $n-2$ elements to 1. The second possibility is if $\alpha = \alpha^{(i)}_{0,p}$ and $\beta =  \alpha^{(j)}_{0,q}$, where $p, q = 1, \ldots, n-2$ and $p \neq q$. Then $\alpha + \beta$ transforms $n-3$ elements to 1. Since $k > 2$, it follows $n-3 \geq n-k$, so, in both cases $\alpha + \beta \in I$.

\vspace{1mm}

\emph{Case 3.} Let $\alpha(a_i) = 0$ for some $i = 1, \ldots, k$ and $\alpha(1) = 1$. Then $\alpha$, represented as an ordered $n-1$ -- tuple, has $n -2$ coordinates equal to 1. For any endomorphism $\beta$ it follows that $n-1$ -- tuple $\alpha + \beta$ has either $n -2$ coordinates equal to 1, or $\alpha + \beta = \overline{1}$. Hence, $\alpha + \beta \in I$.

\vspace{1mm}

\emph{Case 4.} Let  $\alpha(a_i) = 1$ for some $i = 1, \ldots, k$.  Then it follows $\alpha(1) = 1$. From the same arguments, as in Case 3, it follows that $\alpha + \beta \in I$.

Thus we prove that $I$ is closed under the addition.

\vspace{1mm}

Now let $\alpha \in I$ and $\beta \in {\mathcal{E}}_{\mathcal{\lozenge}_n}(A)\backslash I$. Since  $\alpha$ transforms $n-k$ elements or more to 1 and $\beta(1) = 1$, then for any element $x$ of this sort it follows $(\alpha\cdot\beta)(x) = \beta(\alpha(x)) = \beta(1) = 1$. Hence, we have $\alpha\cdot \beta \in I$.

In order to prove that $\beta\cdot \alpha \in I$ we choose $\alpha \in {\mathcal{E}}_{\mathcal{\lozenge}_n}(A_j)$ for some $j = 1, \ldots, k$. Then $a_j \notin \im(\alpha)$. Since $a_j \in \im(\beta)$, there is $a_i$, $i = 1, \ldots, n-2$ such that $\beta(a_i) = a_j$. Now we shall consider three cases.

\vspace{1mm}
\emph{Case 5.} Let $\alpha(a_j) = 0$.  Then $(\beta\cdot \alpha)(a_i) = \alpha(\beta(a_i)) = \alpha(a_j) = 0$ and $\beta\cdot \alpha \in I$.

\vspace{1mm}

\emph{Case 6.} Let $\alpha(a_j) = a_k$, where $k \neq j$. Then $(\beta\cdot \alpha)(a_i) = \alpha(\beta(a_i)) = \alpha(a_j) = a_k$ and $\beta\cdot \alpha \in I$.
\vspace{1mm}

\emph{Case 7.} Let $\!\alpha(a_j) = 1$. We know that  endomorphism $\beta$ transforms $n-k-1$ or more ele\-ments to 1. Then, since $\alpha(1) = 1$, it follows that $\beta\cdot \alpha$ also transform all these $n-k-1$ elements to 1. But $\beta(a_i) = a_j$ and then $(\beta\cdot \alpha)(a_i) = \alpha(\beta(a_i)) = \alpha(a_j) = 1$. So, $\beta\cdot \alpha$  transform $n - k$ elements to 1. Hence, $\beta\cdot \alpha \in I$.

\vspace{1mm}

 Finally, let us observe that the elements of the set ${\mathcal{E}}_{\mathcal{\lozenge}_n}(A)\backslash I$ are isomorphic to permutations of the  $A$. This implies that $I$ is a maximal ideal of ${\mathcal{E}}_{\mathcal{\lozenge}_n}(A)$.

\vspace{3mm}

Note that Proposition 4.4 is a particular case of the last theorem for $k = n-2$.

\vspace{7mm}

{\bf \hspace{6.4mm} References}

\vspace{3mm}

[1] A. Anderson and N. Belnap, \emph{ Entailment, the Logic of Relevance
and Necessity}, vol. I, Princeton Univ. Press, Princeton, 1975.

[2] R. El Bashir and T. Kepka, \emph{Congruence-Simple Semirings}, Semigroup Forum, Vol. 75 (2007) 588 -- 608.

[3]  J. Golan, \emph{Semirings and Their Applications}, Kluwer, Dordrecht, 1999.

[4] G. Gratzer, \emph{Lattice Theory: Foundation}, Birkh\"{a}user
Springer Basel AG, 2011.

[5]  J. Je$\hat{\mbox{z}}$ek, T. Kepka and M. Mar\`{o}ti, \emph{The endomorphism semiring of a se\-milattice},
Semigroup Forum, 78 (2009), 21 -- 26.

[6]  G. Maze, C. Monico and J. Rosenthal, \emph{A public key cryptosystem based
on actions by semigroups}, Advances in Mathematics of Communications,
Volume 1, No. 4, 2007, 489 - 507.

[7] C. Monico, \emph{On finite congruence-simple semirings}, J. Algebra 271 (2004)
846 -- 854.

[8]  I. Trendafilov and  D. Vladeva, \emph{The endomorphism semiring of a finite chain}, Proc. Techn. Univ.-Sofia, 61, 1, (2011), 9 -- 18.

[9]  I. Trendafilov and  D. Vladeva, \emph{Subsemirings of the endomorphism semiring of a finite chain}, Proc. Techn. Univ.-Sofia, 61, 1, (2011), 19 -- 28.

[10]  I. Trendafilov and D. Vladeva, \emph{Endomorphism semirings without zero of a finite semilattice of a special type}, Proc. Techn. Univ.-Sofia, 61, 2, (2011), 19 -- 28.

[11] I. Trendafilov and  D. Vladeva, \emph{Idempotent elements of the endomorphism semiring of a finite chain}, ISRN Algebra, Volume 2013 (2013) (to appear)

[12] J.  Zumbr\"{a}gel, \emph{Classification of finite congruence-simple semirings with zero},
J. Algebra Appl. 7 (2008) 363 -- 377.

\vspace{10mm}

Department of Algebra and Geometry, Faculty of Applied Mathematics and Informatics,\break Technical University of Sofia, 8 Kliment Ohridski Str. Sofia 1000, Bulgaria

\vspace{2mm}

 \emph{e-mail:} ivan$\_$d$\_$trendafilov@abv.bg

\end{document}